\documentclass{elsart}
\usepackage{natbib}
\usepackage{amsmath}
\usepackage{amssymb}
\usepackage{graphics}
\usepackage{graphicx}
\usepackage{multirow}
\usepackage{color}

\def \R{\mathrm{I}\!\mathrm{R}}

\begin{document}

\begin{frontmatter}
	\title{Nonlinear Simplex Regression Models}
	\author[1]{Patr\'{\i}cia L.\ Espinheira\corauthref{cor}},
	\corauth[cor]{Corresponding author.}
	\ead{patespipa@de.ufpe.br}
	\author[1]{Alisson de Oliveira Silva}
	\address[1]{Departamento de Estat\'{\i}stica, Universidade Federal de Pernambuco,
		Cidade Universit\'aria, Recife/PE, 50740--540, Brazil}
\begin{abstract}
\normalsize
\noindent 
In this paper, we propose a simplex regression model in which both the mean and the dispersion parameters are related to covariates by nonlinear predictors. We provide closed-form expressions for the score function, for Fisher's information matrix and its inverse. Some diagnostic measures are introduced. We propose a residual, obtained using Fisher's scoring iterative scheme for the estimation of the parameters that index the regression nonlinear predictor to the mean response and numerically evaluate its behaviour. We also derive the appropriate matrices for assessing local influence on the parameter estimates under diferent perturbation schemes. We also proposed a scheme for the choice of starting values for the Fisher's  iterative scheme for nonlinear simplex models.
The diagnostic techniques  were applied on actual data. The local influence analyses reveal that the simplex models can be a modeling alternative more robust to influential cases than the beta regression models, both to linear and nonlinear models.
\end{abstract}


\end{frontmatter}

\section{Introduction}
In several practical situations, whether experimental or observational, there is interest in investigating how a set of variables relates to percentages or rates.  Among the most suitable models for these data types are the simplex and beta regression models \citep{KieschnickMcCullough2003}.
The beta regression model was proposed by \cite{FERRARI2004} and extended to nonlinear models by \cite{SIMAS2010}.  
\cite{ROCHA2011} also developed residual and local influence for the class of nonlinear beta regression models.

A competitive alternative to the beta regression model is the  simplex regression model proposed by \cite{Barndorff}. The simplex distribution is part of the dispersion models defined by \cite{Jorgensen1997} that extends the generalized linear models. The simplex distribution has been widely used to model data constrained to the interval $(0,1)$. \cite{SONG2000}, for example, proposed a simplex regression model with constant dispersion to model longitudinal continuous proportion data under the generalized estimation equation (EEG) approach. This approach was modified by \cite{SONG2004}, assuming that the dispersion parameter varying throughout the observations. After that, 
\cite{QIUSONGTAN2008}  introduced the
simplex mixed-effects model and more recently \cite{PengQiuShi2016} implemented the package simplexreg in the R system  and available from the Comprehensive R Archive Network (CRAN) at https:
//CRAN.R-project.org/package=simplexreg.  

In this paper, as our most important contribuition we propose the general nonlinear simplex regression model, which considers nonlinear structures in the parameters, for both the mean and dispersion submodels. 
We present analytical expressions for the score vector, for the Fisher information matrix and its inverse, and the estimation of the model parameters is performed by maximum likelihood.  To that end, arises the second important contribuition of this paper, the proposal of a scheme for choosing initial guesses for the maximum likelihood estimation of the parameters under nonlinearity for the simplex regression model. Our proposal is based on \cite{ESPINHEIRA2017}, which has been shown to be extremely relevant on ensuring the convergence of the log-likelihood maximization  procedure  
based on  BFGS method with analytical first derivatives for nonlinear beta regression models.

We also propose diagnostic tools for the nonlinear simplex regression model, being our third  important contribution.  We have developed a residual based on Fisher's iterative scoring algorithm for
estimating of the  mean submodel coefficients.
Sometimes, the distribution of this residual is not well approximated by standard normal distribuition. Thus, we decided to to adopt strategy of using
new thresholds for residual plots based on the simulated envelope algorithm which are associated with normal probability plots proposed by \cite{ESPINHEIRA2017}. The authors show that by using these thresholds, one is more likely to identify atypical points than when using limits ($-2$ and $2$) based on the normal distribution assumption for the residuals.

Finally, we developed local influence (Cook~1986) based on two traditional pertubations schemes, known as: case
weighting and covariate perturbation. Furthermore, 
we propose a new perturbation approach for the response variable, that considered the peculiar aspects of the simplex distribution. Those local influence schemes are appropriate both for linear and nonlinear simplex regression models.

We presented two applications to real data in order to evaluate the behaviour of simplex regression models
relative to the beta regressions models. We show that both models are important and one can be an alternative to another depending on the features of the response variable. In particular, in the two examples presented, the modeling based on the simplex distribution presented a process of estimation by maximum likelihood more robust to influential points than the models based on the beta distribution.
\section{Nonlinear beta regressions}\label{L:model} 
If a random variable $y$ follows the simplex distribution, denoted by $\mathcal{S}^{-1}(\mu, \sigma^2)$ with parameters $\mu \in (0,1)$ and $\sigma^2 > 0$, the density is given by
\begin{equation}\label{denssimplex}
p(y;\mu,\sigma^2) = [2\pi\sigma^2\{y(1-y)\}^3]^{-1/2}\exp\left\{-\frac{1}{2\sigma^2}d(y;\mu)\right\}, y \in (0,1),
\end{equation}  
\noindent in which the  deviation component $d(y;\mu)$ is given by 
\begin{equation}\label{deviance}
d(y;\mu) = \frac{(y-\mu)^2}{y(1-y)\mu^2(1-\mu)^2}.
\end{equation}

The variance function for the simplex distribution is given by
$V(\mu) = \mu^3(1-\mu)^3$. By \cite{Jorgensen1997}, it follows that $\mathrm{E}(y) = \mu$ and 
\begin{equation}\label{VarResposta}
\mathrm{Var}(y) = \mu(1-\mu) - \sqrt{\frac{1}{2\sigma^2}}\exp\left\{\frac{1}{2\sigma^2\mu^2(1-\mu)^2}\right\}\Gamma\left\{\frac{1}{2},\frac{1}{2\mu^2\sigma^2(1-\mu)^2}\right\}.
\end{equation}
\noindent Here, $\Gamma(a,b)$ is the incomplete gamma function defined by
$\Gamma(a,b) = \int_{b}^{\infty}t^{a-1}\,\mathrm{e}^{-t}dt$.

The simplex distribution is very flexible to model data in the continuous interval $(0,1)$, presenting different shapes according to the values of the parameters that index the distribution. In addition to the asymmetric left, asymmetric right, $J$, $U$ 
and $J$ inverted forms, known from the beta distribution, the simplex model is very useful to accommodate data with bimodal distributions. 
\section{The class of nonlinear simplex regression model}\label{R:SimplexNL}

\noindent Let $y_1, \ldots, y_n$ be independent random variables such that 
each $y_t$, $t=1,\ldots,n$, is simplex-distributed, i.e., each $y_t$ has density (\ref{denssimplex})
where $0 < \mu_t < 1$ and $\sigma^2_t > 0$. The nonlinear simplex regression models are defined by (\ref{denssimplex}) 
and by the systematic components:
\begin{eqnarray}\label{compsistem}
g(\mu_t)=f_1(x_t^{\top};\beta) = \eta_t \quad \mbox{and} \quad h(\sigma^{2}_{t}) = f_2(z_t^{\top};\gamma) = \zeta_t,
\end{eqnarray}
\noindent in which $\beta = (\beta_1, \ldots, \beta_k)^{\top}$ and $\gamma = (\gamma_1, \ldots, \gamma_q)^{\top}$ are, respectively,   
a $k$-vector and a $q$-vector of unknown parameters $\beta \in \R^{k}$ and $\gamma \in \R^q$, $k + q < n$, $\eta = (\eta_1,\ldots,\eta_n)^{\top}$ and $\zeta = (\zeta_1,\ldots,\zeta_n)^{\top}$ are nonlinear predictors, $f_1(\cdot)$ and $f_2(\cdot)$ are linear or nonlinear functions, continuous and differentiable in the second argument, so that the matrices of derivatives $\tilde{\mathcal{X}} = \partial \eta/ \partial \beta$ and $\tilde{\mathcal{Z}} = \partial \zeta/ \partial \gamma$  have ranks $k$ and $q$, respectively.
Here, for $t=1,\ldots,n$, $x_t^{\top}=(x_{t1}, \ldots, x_{tk_1})$ and  $z_t^{\top}=(z_{t1}, \ldots, z_{tq_1})$ are, respectively, $k_1$ and $q_1$ observations of known covariates, which may coincide in whole or in part, such that  $k_1\leq k$ and $q_1\leq q$. We assume that the link functions $g:(0,1)\to \R $ and $ h:(0, \infty)\to \R$ are strictly monotonous and twice differentiable.
Different link functions can be chosen for being $g$ and $h$. For example, for $\mu$ we can use the logit specification: $g(\mu) = \log\{\mu/ (1- \mu)$; the probit function: $g(\mu) = \Phi(\mu)$, where $\Phi(\cdot)$ denotes the standard normal distribution function;  the complementary log-log function  $g(\mu) = \log\{-\log(1-\mu)\}$ and the log-log function $g(\mu) = \log\{-\log(\mu)\}$, among others. Since that $\sigma^2 > 0$, one can use $h(\sigma^2) = \log(\sigma^2)$ or the square root function  $h(\sigma^2) = \sqrt{\sigma^2}$. Even the  identity function $h(\sigma^2) = \sigma^2$ can be used, since that one takes into account the positiveness of the estimatives. 

Based on \ref{denssimplex}, it follows that the logarithm of the likelihood function has the form 
$\ell(\beta,\gamma) = \sum_{t = 1}^{n}\ell_{t}(\mu_t,\sigma^{2}_{t})$, 
\noindent in which

\begin{equation}\label{logverossimilhanca}
\ell_t(\mu_t,\sigma^{2}_{t}) = -\frac{1}{2}\log 2\pi - \frac{1}{2}\log \sigma^{2}_{t} - \frac{3}{2}\log\{y_t(1-y_t)\} - \frac{1}{2\sigma^{2}_{t}}d(y_t;\mu_t).
\end{equation}

Thus, follows that the score vector $(U_{\beta}(\beta,\gamma)^{\top},U_{\gamma}(\beta,\gamma)^{\top})^{\top}$ is defined  by
\begin{equation}\label{Escore_beta}
U_{\beta}(\beta,\gamma)= \tilde{\mathcal{X}}^{\top}SUT(y-\mu) \quad
\text{and}\quad
U_{\gamma}(\beta,\gamma)=\tilde{\mathcal{Z}}^{\top}Ha.
\end{equation}
In \eqref{Escore_beta}, $y = (y_1,\ldots,y_n)^{\top}$, $\mu = (\mu_1,\ldots,\mu_n)^{\top}$, $a = (a_1,\ldots,a_n)^{\top}$ are $n$ vectors and $U = \mathrm{diag}\{u_1, \ldots, u_n\}$, where 
\begin{equation}\label{compveca}
a_t = \left\{\frac{d(y_t;\mu_t)}{2(\sigma_{t}^{2})^2}-\frac{1}{2\sigma_{t}^{2}}\right\} \quad \mbox{and} \quad
u_t = \frac{1}{\mu_t(1-\mu_t)}\left\{d(y_t;\mu_t) + \frac{1}{\mu_t^2(1-\mu_t)^2}\right\}. 
\end{equation}
with $d(y_t;\mu_t)$ given by \eqref{deviance}.
Furthermore, 
\begin{equation}\label{Smatrix}
S = \mathrm{diag}\left\{\frac{1}{\sigma^2_1}, \ldots, \frac{1}{\sigma^2_n}\right\},\quad
T = \mathrm{diag}\left\{\frac{1}{g'(\mu_1)}, \ldots, \frac{1}{g'(\mu_n)}\right\}\quad \text{and} \quad
\end{equation}
$H = \mathrm{diag}\left\{\frac{1}{h'(\sigma^2_1)}, \ldots, \frac{1}{h'(\sigma^2_n)}\right\}.$
Fisher's information for $\beta$ and $\gamma$ $K(\beta,\gamma)$ is a  diagonal block matrix, in which
$K_{\beta\beta} = \tilde{\mathcal{X}}^{\top}S W\tilde{\mathcal{X}}$ and $K_{\gamma\gamma} = \tilde{\mathcal{Z}}^{\top}D\tilde{\mathcal{Z}}$. Here,  $W = \mathrm{diag}\{w_1,\ldots,w_n\}$ and $D = \mathrm{diag}\{d_1,\ldots,d_n\}$, with
\begin{equation}\label{compW}
w_t= \frac{ v_t }{\sigma^2_tg^{\prime}(\mu_t)^2},\, v_t= \sigma_t^{2}\left\{\frac{3\sigma_t^{2}}{\mu_t(1 - \mu_t)} + \frac{1}{\mu_t^3(1 - \mu_t)^3}\right\} \, \mbox{and} \,\,
d_t = \left\{\frac{1}{2(\sigma_{t}^{2})^2 h^{\prime}(\sigma_{t}^{2})^2}\right\}.
\end{equation}
Since $K(\beta,\gamma)$ is a diagonal block matrix, the vectors $\beta$ and $\gamma$ are globally orthogonal \citep{CoxReid1987}, so that their maximum likelihood estimators $\widehat{\beta}$ and $\widehat{\gamma}$, respectively, are asymptotically independent. For large samples and under regularity conditions,  we have approximately that
\begin{eqnarray*}
	\left(
	\begin{array}{c}
		\widehat{\beta}\\
		\widehat{\gamma}
	\end{array}
	\right) \sim 
	N_{k + q}
	\left(
	\left(
	\begin{array}{c}
		\beta \\
		\gamma
	\end{array}
	\right),K^{-1}
	\right),\quad {\mbox{with}}\quad K^{-1} = K(\beta,\gamma)^{-1} =
	\left( 
	\begin{array}{cc} 
		K^{\beta\beta} & 0 \\
		0 & K^{\gamma\gamma}
	\end{array}
	\right).
\end{eqnarray*}
Here, 
\begin{equation}\label{invFisher}
K^{\beta\beta} = (\tilde{\mathcal{X}}^{\top} SW\tilde{\mathcal{X}})^{-1}\quad \mbox{and} \quad K^{\gamma\gamma} = (\tilde{\mathcal{Z}}^{\top}D\tilde{\mathcal{Z}})^{-1}.
\end{equation}
The maximum likelihood estimators of $\beta$ and $\gamma$ are obtained as the solution of  $U_{\beta}(\beta, \gamma) = 0$ and $U_{\gamma}(\beta,\gamma) = 0$. However, we emphasize that the maximum likelihood estimators in nonlinear models do not usually present closed-form analytical expressions, making it necessary to use iterative methods, such as quasi-Newton algorithms (eg BFGS); See for example \cite{NocedalWright1999}. The optimization algorithms require the specification of a value $\theta^{(0)} = ({\beta^{(0)}}^{\top}, {\gamma^{(0)}}^{\top})^{\top}$ in order to start the iterative process.
\section{New starting values for log-likelihood maximization}\label{S:startingpoint}
The estimation of parameters in nonlinear models can be an arduous task when Fisher's scoring method is used. However, a well-structured starting-values process, taking into account the features of the response variable distribution, can guarantee both the convergence of the process and the feasibility of the maximum likelihood estimates.
Our proposal is based on \cite{ESPINHEIRA2017} that proposed a starting-values procedure for the nonlinear beta regression model.
We shall consider that $k_1 = k$ and $q_1 = q$. Thus, 
we take the first order Taylor expansion of $f(x_t,\beta)$ at $\beta^{(0)}$, given by
\begin{eqnarray}
f(x_t,\beta) \approx f(x_t,\beta^{(0)}) + \sum_{t = 1}^{k}\left[\frac{\partial f(x_t,\beta)}{\partial \beta_t}\right]_{\beta = \beta^{(0)}}(\beta_t - \beta_t^{(0)}),
\end{eqnarray}
\noindent where $\beta^{(0)} = (\beta_{1}^{(0)},\ldots,\beta_{k}^{(0)})$ is an initial guess. Thus, $f(x_t,\beta) = f(x_t,\beta^{(0)}) + \sum_{i = 1}^{k}\tilde{x}_{ti}^{(0)}(\beta_i - \beta_i^{(0)})$. 
We consider $f(x_t,\beta) = g(y_t)$,  $\theta_i^{(0)} = (\beta_i - \beta_i^{(0)})$, thus $g(y_t) - f(x_t,\beta^{(0)}) = \sum_{i = 1}^{k}\tilde{x}_{ti}^{(0)}\theta_{i}^{(0)}$. With the perspective of a linear model we have that the least squares estimator of $\theta^{(0)}$ is given by $\widehat{\theta}^{(0)} = (\tilde{\mathcal{X}}^{(0)\top}\tilde{\mathcal{X}}^{(0)})^{-1}\tilde{\mathcal{X}}^{(0)\top}(g(y)-f(x,\beta^{(0)}))$, where $\tilde{\mathcal{X}}^{(0)} = [\partial \eta/\partial \beta]_{\beta = \beta^{(0)}}$ and $\widehat{\theta}_{i}^{(0)} = (\widehat{\beta}_{i}^{(1)} - \beta_{i}^{(0)})$. Hence, $\widehat{\beta}_{i}^{(1)} = \widehat{\theta}_{i}^{(0)} + \beta_{i}^{(0)}$. Our  proposal is to use is use the following nonlinear initial guess to $\widehat\beta$: $\beta_{NL}^{(0)} = (\tilde{\mathcal{X}}^{(0)\top}\tilde{\mathcal{X}}^{(0)})^{-1}\tilde{\mathcal{X}}^{(0)\top}(g(y)-f(x,\beta_{L}^{(0)}))$ and $\beta^{(0)} = \beta_{L}^{(0)} = (X^{\top}X)^{-1}X^{\top}g(y)$. As can be noticed, the procedure we propose for
the starting points involves two steps. The first step is obtain $\beta_{L}^{(0)}$ and then, we can compute the nonlinear initial guess $\beta_{NL}^{(0)}$.
For the precision/dispersion submodel, we consider $h(\sigma_t^2) = f(z_t,\gamma)$. Hence,
$\gamma_{NL}^{(0)} = (\tilde{\mathcal{Z}}^{(0)\top}\tilde{\mathcal{Z}}^{(0)})^{-1}\tilde{\mathcal{Z}}^{(0)\top}(h(\sigma_{NL}^{2})-f(z,\gamma_{L}^{(0)}))$ and $\gamma^{(0)} = \gamma_{L}^{(0)} = (Z^{\top}Z)^{-1}Z^{\top}h(\sigma_{L}^{2(0)})$.
Here, $\tilde{\mathcal{Z}}^{(0)} = [\partial \zeta/\partial \gamma]_{\gamma = \gamma^{(0)}}$ and $\sigma_{tL}^2 = d(y_t;\check{\mu}_{Lt})$ with $\check{\mu}_{Lt} = g^{-1}(\hat{\eta}_{1L}) = g^{-1}(x_{t}^{\top}\beta_{L}^{(0)})$. Finally, $\sigma_{NLt}^{2(0)} = d(y_t;\check{\mu}_{NLt})$ with $\check{\mu}_{NLt} = g^{-1}(\hat{\eta}_{1NL}) = g^{-1}(f(x_t^{\top},\beta_{NL}^{(0)}))$.

Typically, in one or both submodels there are more parameters than covariates in nonlinear predictors, that is,  $k_1< k$ and/or $q_1 < q$. In those cases, it is necessary given  numerical values for parameters to obtain a predictor formed by covariates that no longer involve unknown parameters. This step is very relevant, since that these initial numerical guess shall respect the features for both the covariates and its relation with the mathematical functions that describe the nonlinear predictors. After this step we shall construct a matrix $X$ based on a linear predictor and compute $(X^{\top}X)^{-1}X^{\top}g(y)$. Based on this two step procedure it is possible to obtain $\beta_{L}^{(0)}$ when numbers of parameters exceed the number of covariates.
\section{Weighted residual}\label{S:Weightedresidual} 
We can define residual as a measure that aims to identify discrepancies between the estimated model and the data. Thus, it is relevant to take in accout the probability distribution's features of the response variable, as well as the features of the estimation process of the model. From this perspective, \cite{ESPINHEIRA2008Res} suggested to use standardized residuals obtained from the convergence of the iterative Fisher's scoring process for estimation of the vector of the regression parameters.

Based on the nonlinear simplex regression model in \eqref{compsistem}, on the score function for $\beta$ in \eqref{Escore_beta} and on the inverse of Fisher's information for $\beta$ in \eqref{invFisher}, it follows that Fisher's scoring iterative algorithm for 
estimating $\beta$ is given by 
\begin{eqnarray}\label{eq:iterativo}
\beta^{(m + 1)} = \beta^{(m)} + (\tilde{\mathcal{X}}^{\top(m)}S^{(m)} W^{(m)}\tilde{\mathcal{X}}^{(m)})^{-1}\tilde{\mathcal{X}}^{\top(m)}S^{(m)}U^{(m)}T^{(m)}(y-\mu^{(m)}), \quad\,\,
\end{eqnarray}
where $m=0,1,2,....$ indexes the iterations that take place until 
convergence, which occurs when the distance between $\beta^{(m+1)}$ and $\beta^{(m)}$ becomes 
smaller than a given small constant. 
It is possible to write the iterative scheme in (\ref{eq:iterativo}) in terms of weighted least squares 
regressions, such that
$\beta^{(m + 1)}=(\tilde{\mathcal{X}}^{\top(m)}S^{(m)} W^{(m)}\tilde{\mathcal{X}}^{(m)})^{-1}\break\tilde{\mathcal{X}}^{\top(m)}S^{(m)}W^{(m)}z^{(m)}$, where $z^{(m)} = \tilde{\mathcal{X}}^{(m)}\beta^{(m)} + W^{-1(m)}U^{(m)}T^{(m)}(y-\mu^{(m)})$. Upon convergence,
\begin{eqnarray}
\widehat{\beta} = (\widehat{\tilde{\mathcal{X}}}^{\top}\widehat{S}\widehat{W}\widehat{\tilde{\mathcal{X}}})^{-1}\widehat{\tilde{\mathcal{X}}}^{\top}\widehat{S}\widehat{W}z,\quad \text{with}\quad	z = \widehat{\tilde{\mathcal{X}}}\widehat{\beta} + \widehat{W}^{-1}\widehat{U}\widehat{T}(y-\widehat{\mu}).\label{betachapeu}
\end{eqnarray}
We have that $\widehat{\beta}$ in (\ref{betachapeu}) can be accounted as the least squares
estimate of $\beta$ obtained by regressing $z$ on $\widehat{\tilde{X}}$ with weighting matrix $\widehat{S}\widehat{W}$. 
In this sense, the ordinary residual is 
\begin{align} r^{*} &= (\widehat{S}\widehat{W})^{1/2}(z - \widehat{\tilde{X}}\widehat{\beta}) = [(\widehat{S}\widehat{W})^{1/2} - (\widehat{S}\widehat{W})^{1/2}\widehat{\tilde{X}}(\widehat{\tilde{\mathcal{X}}}^{\top}\widehat{S}\widehat{W}\widehat{\tilde{\mathcal{X}}})^{-1}\widehat{\tilde{\mathcal{X}}}^{\top}\widehat{S}\widehat{W}]z\label{residuo} \\ 
&=\widehat{S}^{1/2}\widehat{W}^{-1/2}\widehat{U}\widehat{T}(y-\widehat{\mu})\nonumber
\end{align}
An alternative for standardize the ordinary residual is using an aproximation to the variance 
of $z$. Let (\ref{betachapeu}) as 
$(\widehat{\tilde{X}}^{\top}\widehat{S}\widehat{W}\widehat{\tilde{X}})\widehat{\beta} = \widehat{\tilde{X}}^{\top}\widehat{S}\widehat{W}z
$ and  consider that $\widehat{S}\approx S$, $\widehat{W}\approx W$ and $\widehat{\tilde{X}} \approx {\tilde{X}}$. Then, using the fact that $\mathrm{cov}(\widehat{\beta}) \approx (\tilde{X}^{\top}SW\tilde{X})^{-1}$, it follows that $\mathrm{cov}(z) \approx (SW)^{-1}$. 
Thus, based on (\ref{residuo}), we have that
\begin{align*}
&\mbox{Cov}(r^{*}) \approx [({S}{W})^{1/2} - ({S}{W})^{1/2}{\tilde{X}}({\tilde{\mathcal{X}}}^{\top}{S}{W}{\tilde{\mathcal{X}}})^{-1}{\tilde{\mathcal{X}}}^{\top}{S}{W}] \times\mbox{Cov}(z) \times\\
&[({S}{W})^{1/2} - ({S}{W})^{1/2}{\tilde{X}}({\tilde{\mathcal{X}}}^{\top}{S}{W}{\tilde{\mathcal{X}}})^{-1}{\tilde{\mathcal{X}}}^{\top}{S}{W}]^\top\\
&\approx (I- {S}^{1/2}{W}^{1/2}{\tilde{X}}{\tilde{\mathcal{X}}}^{\top}{S}{W}{\tilde{\mathcal{X}}})^{-1}{\tilde{\mathcal{X}}}^{\top}{W}^{1/2}{S}^{1/2}) \approx\mathcal{I}_n - H^{*},
\end{align*}
with $H^{*} =  \widehat{S}^{1/2}\widehat{W}^{1/2}\widehat{\tilde{X}}(\widehat{\tilde{\mathcal{X}}}^{\top}\widehat{S}\widehat{W}\widehat{\tilde{\mathcal{X}}})^{-1}\widehat{\tilde{\mathcal{X}}}^{\top}\widehat{W}^{1/2}\widehat{S}^{1/2}$.

Then,
we proposed the {\bf weighted residual} for nonlinear simplex regression models, which includes the class of linear simplex model, and to be defined as
$r^{\beta} =\widehat{S}^{1/2}\widehat{W}^{-1/2}\widehat{U}\widehat{T}(y-\widehat{\mu})(I - H^*).^{-1/2}$
Finally, considering the definitions of $S$, $W$, $U$, $T$ given 
along (\ref{compveca}) and (\ref{compW}), we have that 	
\begin{equation}\label{weigheted}
r_t^{\beta} =  \frac{\widehat{u}_t(y_t - \widehat{\mu_t})}{\sqrt{\widehat{v}_t(1 - h_{tt}^{*})}},
\end{equation}
where $h^*_{tt}$ denotes the $t$th diagonal element of $H^*$ and $v_t$ is given in (\ref{compW}).
\section{Local influence}\label{S:influence}
Let $\theta ={({\beta}^{\!\top},{\gamma}^{\!\top}})^{\!\top}$ is a $(s=k+q)\times 1$ vector of unknown parameters and $\ell(\theta)$ denotes its log-likelihood function. Now, we introduce a perturbation in the
assumed model through a vector $\delta$, $n\times 1$. Thus, 
$\ell_{\delta}(\theta)$ denote the log-likelihood function of the
perturbed model for a given $\delta$.
The maximum likelihood estimators of $\theta$ for the assumed and
perturbed models are denoted by, respectively,  ${\,\widehat{\!\theta}}$ and ${\,\widehat{\!\theta}}_{\delta}$.
The likelihood displacement
$LD_{\delta} = 2\left\{  \ell(\,\widehat{\!\theta})-
\ell(\,\widehat{\!\theta}_{\delta})\right\}$ can be used to attain the influence
of the perturbation on the maximum likelihood estimate. 
If minor modifications at the postulated model  lead to 
considerable
changes in inferential results, then  further investigations must be carried out about a better model to fit  the data.

Cook~(1986) proposed to analyze the local behaviour of 
$LD_{\delta}$
around $\delta_0$, which represents no perturbation, such that 
$LD_{\delta_0}=0$. The author suggested to evaluate  the curvature of the plot of  
$LD_{\delta_0 + aI}$ against $a$, where $a\in I\!\! R$, $I$ is a unit norm direction.  The interest lies in to find
the direction $I_{\rm max}$ corresponding to the largest curvature 
$C_{\rm max}$.
Observations that are jointly influential can be singled out by the index plot of $I_{\rm max}$. 
Cook~(1986) showed that the normal curvature at direction $I$ is   
$C_I(\theta)=2\vert 
I^{\!\top}{\Delta}^{\!\top}\ddot{\ell\,}^{-1}{\Delta}I\vert$, 
where $\ddot{\ell\,} = {\partial^2 \ell(\,\widehat{\!\theta})/{\partial 
		\theta\partial
		\theta^{\!\top}}}$, and $\Delta$ is an $s\times n$ matrix given by 
$\Delta = {\partial^2 {\ell_{\delta}(\theta)} /{\partial \theta\partial 
		\delta^{\!\top}}}$,
evaluated at $\theta = \,\widehat{\!\theta}$ and $\delta = \delta_0$.
Thus, $C_{\rm max}/2$ is the largest eigenvalue of 
$\bf {-{\Delta}^{\!\top}\ddot{\ell\,}^{-1}{\Delta}}$ and $I_{\rm max}$
is the corresponding eigenvector.

In the other hand, \cite{LesaffreVerbeke1998} proposed a measure for identifying individually influential observations, known as the total local
influence of observation $t$ and defined as 
$C_t=2 \vert {\Delta_t}^{\!\top}\ddot{\ell\,}^{-1}{\Delta_t}\vert$,
where $\Delta_t$ is the $t$th column of $\Delta$. $C_t$  is  
the normal curvature  in the direction of the vector whose $t$th 
component equals one and all other elements are zero. 
Observations such as that $C_t >
2\sum_{t=1}^n C_t/n$ can be taken to be individually influential.
If we part the parameter vector $\theta$ as 
$\theta = (\theta_1^{\!\top}, \theta_2^{\!\top})^{\!\top}$, 
we can access the local influence relative to $\theta_1$, such that 
$C_{I;\theta_1} = \vert I^{\!\top}{\Delta}^{\!\top}(\ddot{\ell\,}^{-1}-
\ddot{\ell\,}_{22}){\Delta}I\vert$ and $C_{t;\theta_1}=2 \vert {\Delta_t}^{\!\top}(\ddot{\ell\,}^{-1}
- \ddot{\ell\,}_{22}){\Delta_t}\vert$,
where
$$ \ddot{\ell\,}_{\theta_2\theta_2} =
{\frac{\partial^2 \ell(\theta)}{ {\partial \theta_2\partial 
			\theta_2^{\!\top}}}}
\quad \mathrm{and} \quad
\ddot{\ell\,}_{22} =
\begin{pmatrix}
0 & 0\cr
0 &{\ddot{\ell\,}_{\theta_2\theta_2}}^{-1}\cr
\end{pmatrix}.$$
In this case, 
$I_{{\rm max};\theta_1}$ is the unit norm eigenvector corresponding 
to the largest eigenvalue of $\bf {-{\Delta}^{\!\top}(\ddot{\ell\,}^{-1}
	-\ddot{\ell\,}_{22}){\Delta}}$.

In what follows, we shall develop local influence measures 
for the class of nonlinear simplex regression
model under three different perturbation schemes, namely: case
weighting,  response perturbation and covariate perturbation. 
To that end we need to obtain $\ddot{\ell\,}$ and $\ddot{\ell\,}^{-1}$. Based on \eqref{logverossimilhanca}, it  
follows that $ \ddot{\ell}_{\beta\beta} = -\tilde{\mathcal{X}}^{\top}SQ\tilde{\mathcal{X}} + [\,\,b_{\beta}^{\top}\,\,][\,\,\tilde{\mathcal{X}}_{\beta}\,\,]$, $ \ddot{\ell}_{\beta\gamma} = - \tilde{\mathcal{X}}^{\top}S^2HTU\mathcal{E}\tilde{\mathcal{Z}}$, $ \ddot{\ell}_{\gamma\beta} = (\ddot{\ell}_{\beta\gamma})^{\top}$ and $ \ddot{\ell}_{\gamma\gamma} = - \tilde{\mathcal{Z}}^{\top}\mathcal{V}\tilde{\mathcal{Z}} + [\,\,b_{\gamma}^{\top}\,\,][\,\,\tilde{\mathcal{Z}}_{\gamma}\,\,]$. Here,
$Q = \mathrm{diag}\{q_1, \ldots, q_n\}$,
\begin{eqnarray}\label{compQ}
q_t = \left\{u_t - (y_t - \mu_t)u_t^{\prime} + (y_t - \mu_t)u_t\frac{g^{\prime\prime}(\mu_t)}{g^{\prime}(\mu_t)}\right\}\frac{1}{\{g^{\prime}(\mu_t)\}^2},
\end{eqnarray}
with
\begin{eqnarray*}
	u_t^{\prime} &=& -\left\{\frac{2(y_t - \mu_t)u_t}{\mu_t(1 - \mu_t)} + \frac{3(1 - 2\mu_t)}{\mu_t^4(1 - \mu_t)^4} + \frac{(1 - 2\mu_t)d(y_t; \mu_t)}{\mu_t^2(1 - \mu_t)^2}\right\}, 
\end{eqnarray*}
$\mathcal{V} = \mathrm{diag}\{\nu_1, \ldots, \nu_n\}$, com 
\begin{eqnarray}\label{compVandE}
\nu_t = d_t + a_t\frac{h^{\prime\prime}(\sigma_t^2)}{\{h^{\prime}(\sigma_t^2)^3}\quad \mbox{and}\quad \mathcal{E} = \mathrm{diag}\{(y_1 - \mu_1), \ldots, (y_n - \mu_n)\}.
\end{eqnarray}
Besides that, $b_{\beta} = STU(y - \mu)$, $\tilde{\mathcal{X}}_{\beta} = (\tilde{\mathcal{X}}_{t})$ is  an array of dimension $n \times k \times k$, being $\tilde{\mathcal{X}}_{t}$ a matrix $k \times k$ with elements given by $\partial^2 \eta_t/\partial \beta_i\partial \beta_p$, $b_{\gamma} = Ha$ and $\tilde{\mathcal{Z}}_{\gamma} = (\tilde{\mathcal{Z}}_{t})$ is an array of dimension  $n \times q \times q$, being  $\tilde{\mathcal{Z}}_{t}$ a matrix  $q\times q$ with elements given by $\partial^2 \zeta_t/\partial \gamma_j \partial \gamma_l$. Finally, $[\,\,\cdot\,\,][\,\,\cdot\,\,]$ represents the  bracket product of a matrix by an array as defined by \cite{Wei1998},~p.188. 

The structure of $\Delta$
for each perturbation scheme is given below. 
For {\bf The cases weighting perturbation} we have that  $\ell_{\delta}(\beta,\gamma) = \sum_{t = 1}^{n}\delta_t\ell_t(\mu_t,\sigma_t^2)$. In this case, $\delta_0 = (1, 1, \ldots, 1)^{\top}$, $\Delta_t  = \partial \ell_t(\theta)/\partial \theta$, and thus
\begin{equation*}
\Delta = 
\left(
\begin{array}{c}
\widehat{\tilde{\mathcal{X}}}^{\top}\widehat{S}\widehat{T}\widehat{U}\widehat{\mathcal{E}}\\
\widehat{\tilde{\mathcal{Z}}}^{\top}\widehat{H}\widehat{\mathcal{A}}
\end{array}
\right),
\end{equation*}
with $\mathcal{A} = \mathrm{diag}\{a_1, \ldots, a_n\}$ being that the quantity $a_t$ and the components of the diagonal matrix $U$  were defined in \eqref{compveca}. Additionaly, the diagonal matrices $S$, $T$, $H$ and ${\mathcal{E}}$ were defined, respectively, in \eqref{Smatrix} and  \eqref{compVandE}.

We now move to  the response perturbation scheme, where  
$y_t(\delta) = y_t + \delta_ts(y_t)$,
with $s(y_t)=\sqrt{ V(\widehat\mu_t)}$, and $V(\widehat\mu_t)=\mu^3_t(1-\mu_t^3)$ as the variance function. {\bf This is a new proposal to the standardization of the scale factor}, given that, in general, we use the standard deviation of $y$. Its worth pointing out that the new proposal presents lesser computational cost because the variance of a random variable with simplex distribution involves computations related to the incomplete gamma function, see \eqref{VarResposta}. In this scheme,  $\delta_0=(0, 0 , \ldots, 0)^{\!\top}$ and
\begin{equation*}
\Delta = 
\left(
\begin{array}{c}
\widehat{\tilde{\mathcal{X}}}^{\top}\widehat{S}\widehat{T}\widehat{M}S_y\\
\widehat{\tilde{\mathcal{Z}}}^{\top}\widehat{S}\widehat{H}\widehat{B}S_y
\end{array}
\right),
\end{equation*}
where $M = \mathrm{diag}\{m_1, \ldots, m_n\}$, $B = \mathrm{diag}\{b_1, \ldots, b_n\}$ and $S_y = \mathrm{diag}\{s(y_1), \ldots, s(y_n)\}$, with
\begin{equation*}
m_t = \frac{1}{y_t(1-y_t)}\left\{\frac{2}{y_t(1-\mu_t)^3} + \frac{(1-3\mu_t)}{\mu_t^2(1-\mu_t)^3} - \frac{1}{2}\frac{\partial d(y_t;\mu_t)}{\partial \mu_t}\right\}\quad \mbox{and}
\end{equation*}
\begin{equation*}
b_t = \frac{1}{2\sigma_t^2y_t(1-y_t)}\left\{d(y_t;\mu_t) + \frac{2(y_t-\mu_t)}{y_t\mu_t(1-\mu_t)^2}\right\}.
\end{equation*}
The third and final perturbation scheme involves the simultaneous perturbation of 
a continuous covariate, say $(x_p^{\top}, z_{p'}^{\top})$, $p=1,\ldots,k_1$ and $p'=1,\ldots,q_1$. We 
replace $x_{tp}$ by $x_{tp}+\delta_t s_{x_p}$ and $z_{tp'}$ by $z_{tp'}+\delta_t s_{z_{p'}}$, where $\delta$ is a 
vector of
small perturbations and $s_{x_p}$ is the standard deviation of $x_{p}$ and 
$s_{z_{p'}}$ is the standard deviation of $z_{p'}$
(Thomas and Cook, 1990). In this scheme,  $\delta_0=(0, 0 , \ldots, 0)^{\!\top}$ and
\begin{equation}\label{DeltaCovariada}
\Delta = 
\left(
\begin{array}{c}
- \widehat{\tilde{\mathcal{X}}}^{\top}\widehat{S}^2\widehat{T}\widehat{H}\widehat{U}\widehat{\mathcal{E}}\widehat{\tilde{Z}}_{\delta} - \widehat{\tilde{\mathcal{X}}}^{\top}\widehat{S}\widehat{Q}\widehat{\tilde{X}}_{\delta} + [\,\,\widehat{b}_{\beta}^{\top}\,\,][\,\,\widehat{\tilde{X}}_{\beta\delta}\,\,]\\
- \widehat{\tilde{\mathcal{Z}}}^{\top}\widehat{S}^2\widehat{T}\widehat{H}\widehat{U}\widehat{\mathcal{E}}\widehat{\tilde{X}}_{\delta} - \widehat{\tilde{\mathcal{Z}}}^{\top}\widehat{\mathcal{V}}\widehat{\tilde{Z}}_{\delta} + [\,\,\widehat{b}_{\gamma}^{\top}\,\,][\,\,\widehat{\tilde{Z}}_{\gamma\delta}\,\,]
\end{array}
\right),
\end{equation}
where $\tilde{\mathcal{X}}_{\delta} = \partial \eta(\delta)/\partial \delta$, $\tilde{\mathcal{Z}}_{\delta} = \partial \zeta(\delta)/\partial \delta$, $\tilde{\mathcal{X}}_{\beta\delta} = (\tilde{\mathcal{X}}_{\delta\,t})$ is an array $n \times k \times n$, being $\tilde{\mathcal{X}}_{\delta\,t}$ a matrix $k \times n$ with elements $\partial^2\eta_t(\delta)/\partial \beta_i \partial \delta_t$, $\tilde{\mathcal{Z}}_{\gamma\delta} = (\tilde{\mathcal{Z}}_{\delta\,t})$ is an array $n \times q \times n$, being   $\tilde{\mathcal{Z}}_{\delta\,t}$   a matrix  $q \times n$ with elements $\partial^2\zeta_t(\delta)/\partial \gamma_j \partial \delta_t$.
In \eqref{DeltaCovariada} the elements of the matrices $Q$ and $\mathcal{V}$ are given by \eqref{compQ} and \eqref{compVandE}, respectively. Other elements as in \eqref{DeltaCovariada} are defined in the paragraph below  expression \eqref{compVandE}.
\section{Simulations}  
\noindent Our aim here is to investigate the distribuition of the {\bf weighted residual} proposed in \eqref{weigheted}  using Monte Carlo experiments. The  experiments were carried out using a nonlinear simplex 
regression model in which 
\begin{eqnarray}\label{ResModel}
\log\frac{\mu_t}{1-\mu_t} = \beta_2 + x_{t2}^{\beta_2} + \beta_3 x_{t3} + \beta_4 x_{t4}; \quad\log(\sigma_t^2) = \gamma_1 + z_{t2}^{\gamma_2},\,t=1,\ldots,n.
\end{eqnarray}
The covariate values were 
obtained as random draws following: $x_{t2} \sim \mathcal{U}(0.5,1.5)$, $x_{t3} \sim \mathcal{U}(0,1)$, $x_{t4} \sim \mathcal{U}(-0.5,0.5)$ and $z_{t2} \sim \mathcal{U}(0.5,1.5)$.  The covariate 
values remained constant throughout the simulations.
We should consider scenarios where the response values
are close to one, scattered  on the standard unit 
interval and close to zero. In this sense, we consider
$\beta = (-2.4,1.4,-1.5,-1.7)$, ($\mu \in (0.02,0.32)$),  $\beta = (-1.7,-1.8,1.2,-1.3)$, ($\mu \in (0.19,0.86)$) and  $\beta = (2.1,-1.5,-1.6,-1.2)$,  ($\mu \in (0.78,0.98)$). 

We measure the intensity of nonconstant dispersion as $\lambda =\{\max\sigma^2_t\}/\{\min\sigma^2_t\}_{t=1,\ldots,n}$ and we consider results to $\lambda \approx 12$, $(\gamma = (-1.3,-1.6))$, $\lambda \approx 45$, $(\gamma = (-1.3,-2.1))$ and $\lambda \approx 128$, $(\gamma = (-1.3,-2.4))$. 
The sample sizes are $n=40,80,120$.  We have to point out that only 
$n=40$ values for the covariates were generated, which were replicated two and three times, respectively, to obtain the other sample sizes. With this, we can guarantee that the intensity of the nonconstant dispersion remains the same for the different sample sizes. The Monte Carlo experiments were carried out based on 10,000 replications.

\begin{figure}[!ht]
	\centering
	\includegraphics[width=0.8\linewidth]{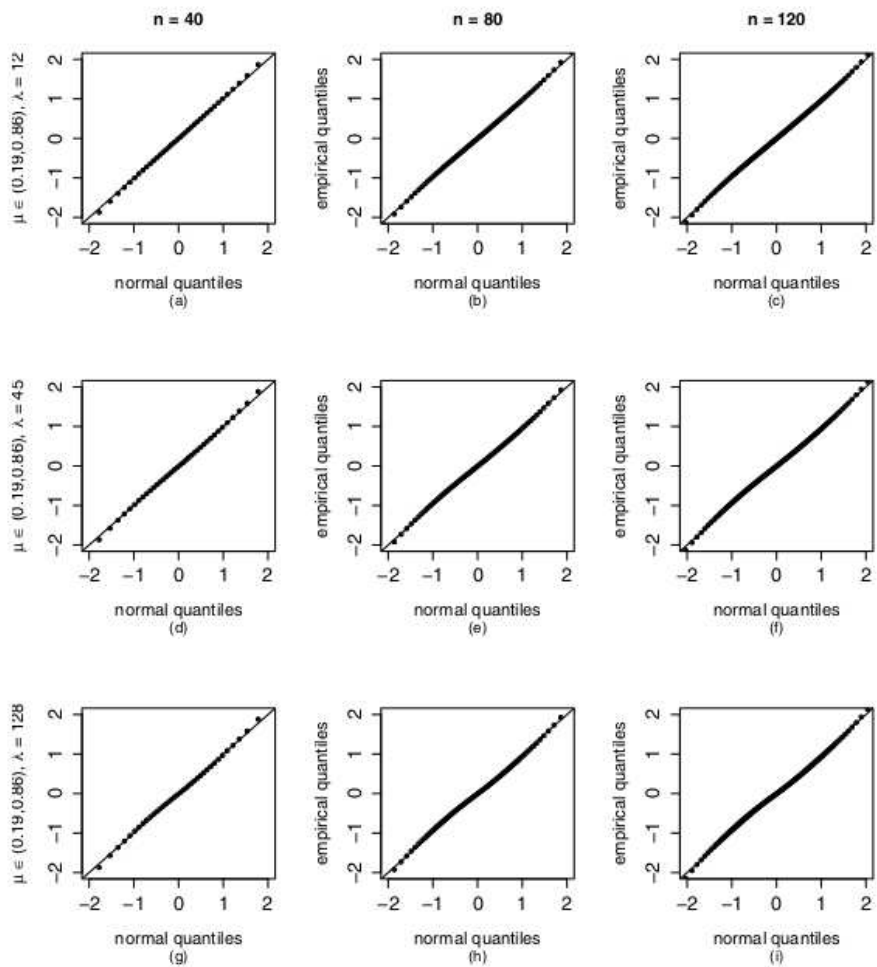}
	\caption{Weighted residuals normal QQ plots. $\mu \in (0.19,0.86)$ }
	\label{fig:QQplotmu2}
\end{figure}

\begin{figure}[!ht]
	\centering
	\includegraphics[width=0.8\linewidth]{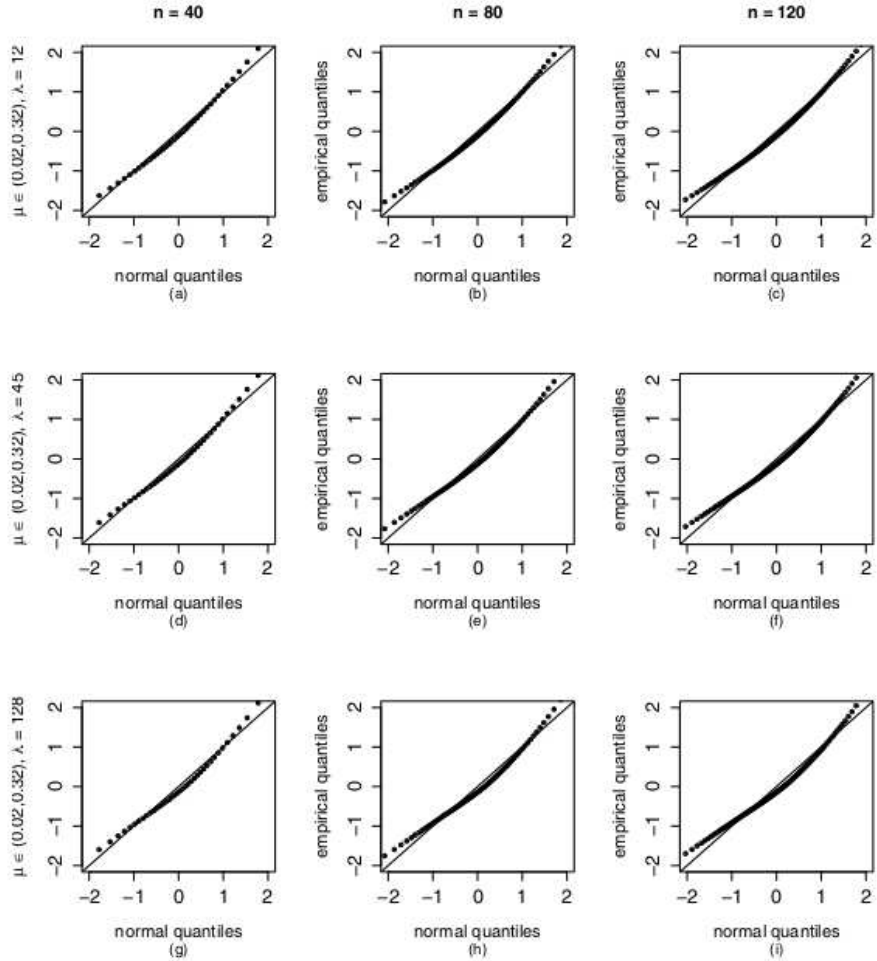}
	\caption{Weighted residuals normal QQ plots. $\mu \in (0.02,0.32)$}
	\label{fig:QQplotmu1}
\end{figure}

\begin{figure}[!ht]
	\centering
	\includegraphics[width=0.8\linewidth]{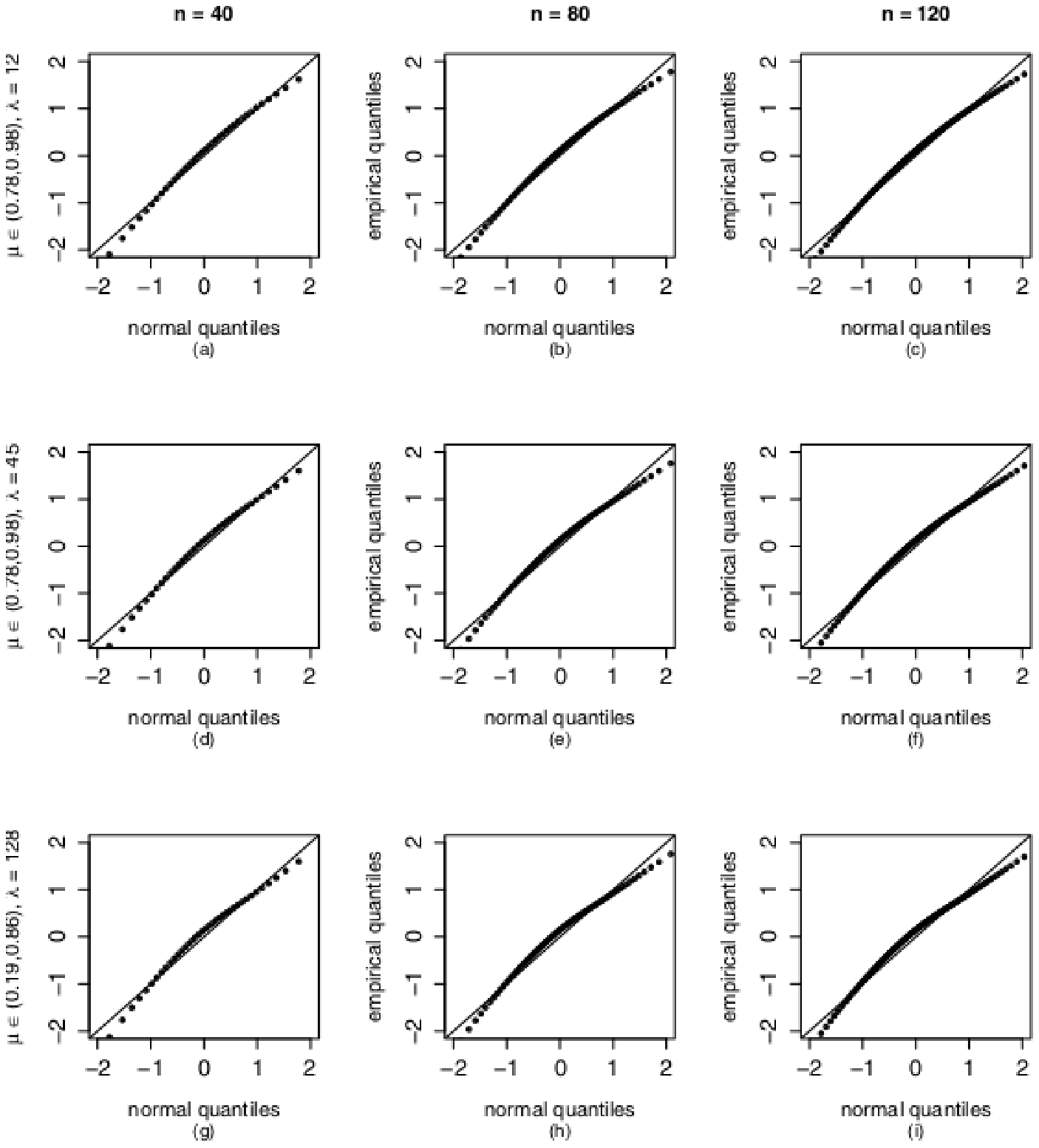}
	\caption{Weighted residuals normal QQ plots. $\mu \in (0.78,0.98)$}
	\label{fig:QQplotmu3}
\end{figure}

The Figures \ref{fig:QQplotmu2}, \ref{fig:QQplotmu1} and \ref{fig:QQplotmu3} contain normal probability plots of the mean order 
statistics of the weighted residual considering the fitted model in \eqref{ResModel} for the different scenarios for $\mu$. 
These figures reveal that the weighted residual distribution is typically well approximated by the standard normal distribution
when
$\mu\in (0.19,0.86)$  even in the case where $\lambda=128$, (Figure \ref{fig:QQplotmu2}). 
However, in cases where $\mu\approx0$ and $\mu\approx 1$, the distribution of the weighted residual presents some asymmetry (Figures \ref{fig:QQplotmu1} and \ref{fig:QQplotmu3}), respectively. Such asymmetry becomes more evident as the intensity of the nonconstant dispersion increases. In addition, for these ranges of $\mu$, the approximation of the distribution of the weighted residual by the standard normal distribution does not seem to improve as the sample size increases.
The acknowledgement of the residual distribution used for the diagnostic analysis is relevant for the definition of detection limits of outliers in the construction of residual plots against the index of the observations, against the predicted values or against covariates. \cite{ESPINHEIRA2017} proposed using empi\-rical  residual quantiles obtained from the simulated envelop bands. Here we used the same strategy.
For any given residual, the thresholds are $\omega_{0.025}$ and $\omega_{0.975}$, the $0.025$ and $0.975$ quantiles of the residual empi\-rical distribution.  
\section{Application I : Reading accuracy data.}
In the first application we consider the data originally analyzed by \cite{SMITHSON2006}. The variable of interest (y) are the scores in a reading accuracy test of 44 children, and the covariates are dyslexia versus 
non-dyslexia status ($x_2$), nonverbal IQ ($x_3$) 
and an interaction variable 
($x_4$). Study participants were recruited from primary schools in the Australian capital. The covariant ($x_2$) assumes the value 1 if the child is dyslexic and -1 otherwise. The mean reading accuracy score was 0.900 for non-dyslexic readers and 0.606 for the dyslexic group. In addition, the scores range from 0.459 to 0.9999, with the overall mean score being 0.773 and the median score being 0.706.

This set of data has been extensively analyzed by several authors. \cite{ESPINHEIRA2008Res}, \cite{ESPINHEIRA2008Inf} and \cite{Ferrari+Espinheira+Cribari_2011} performed an extensive diagnostic analysis to justify the use of a beta regression model in which mean and precision are modeled simultaneously. 
Here we shall consider a simplex regression model with constant dispersion and we  carried out a diagnostic analysis to compare with the diagnostic analysis made by these authors. We consider $g(\mu_t) = \beta_1 + \beta_2 \,x_{t2} + \beta_3 \,x_{t3} + \beta_4 \,x_{t4}$, in which $g(\cdot)$ is the logit function.
\begin{figure}
	\centering
	\includegraphics[width=0.67\linewidth, angle=270]{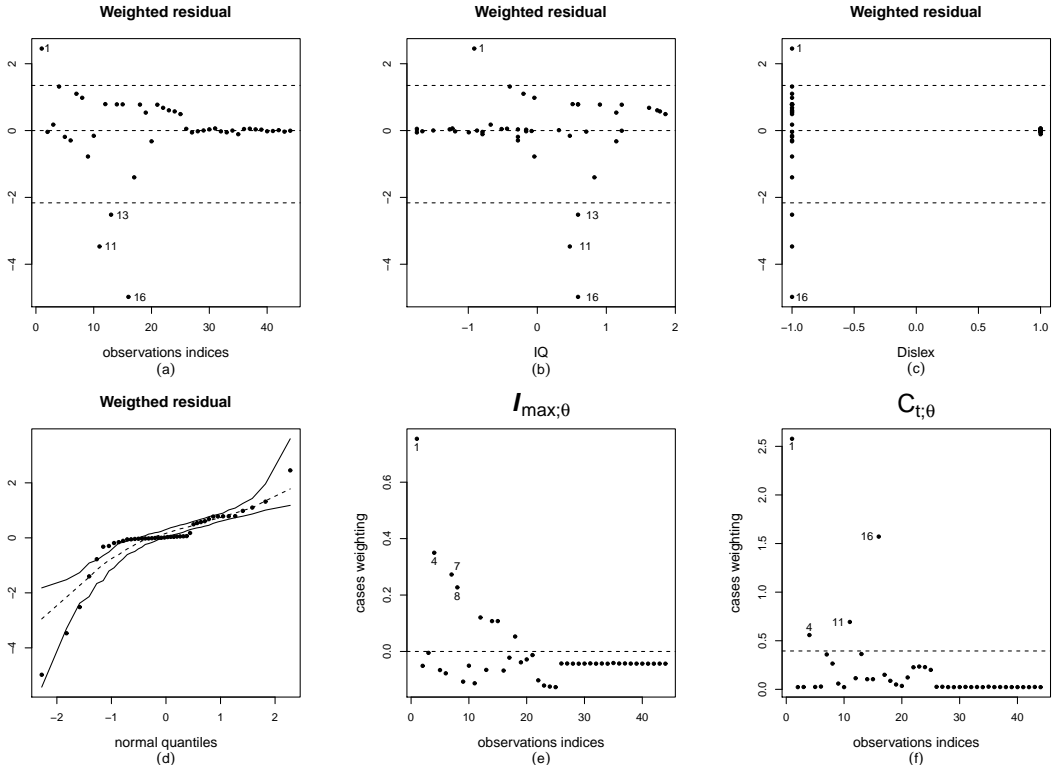}
	\caption{Diagnostic plots and influence plots. Reading accuracy data. Simplex model with constant dispersion.}
	\label{fig:DislexResidual}
\end{figure}
Based on residual and influential plots (Figure \ref{fig:DislexResidual}), we reach the same conclusions found for a beta regression model with constant dispersion fitted for these data by \cite{ESPINHEIRA2008Res}, \cite{ESPINHEIRA2008Inf} and \cite{Ferrari+Espinheira+Cribari_2011}.
Even the same influential observations are singled out. Thus, we present in Table~\ref{tab1}   a comparison between beta and simplex regression models with constant dispersion, considering the full dataset and sub datasets due to exclusions of cases singled out as potentially influential. Furthermore, the need of the modeling of the dispersion is clear. This feature is  more clear by the weighted residual plots of the simplex regression model than what was shown on the residual plots of the beta regression model presented by \cite{ESPINHEIRA2008Res}, (See Figure \ref{fig:DislexResidual} (a)).  
In fact, there is more dispersion in the control group than in  the dyslexic 
group (Figure \ref{fig:DislexResidual}(c)). This fact explains the trend presented in Figure \ref{fig:DislexResidual}(a), since that the first 25 cases are about non-dilex group. Also, IQ covariate seems to be associated with the data dispersion. The normal probability and influence plots also revealed that there is a  problem  associated with the assumption of  constant dispersion.

However, the most important conclusions are presented in Table~\ref{tab1}.
The ana\-lysis of this table reveals that the estimation process of the simplex model is more robust to influential cases than the estimation process of the beta regression model.  It is noteworthy how the case 1 does not affect the inference of the simplex regression model as it affects the inference of the beta regression model. At simplex regression model, even without modeling of the dispersion, both the IQ and the interaction term are covariates with statistically significant effects to explain the mean response. 

We have to stress the feature os this dataset, in which there is an evidence 
of an observation  highly influential. 
Note that the estimate 
of the precision parameter is considerably small $\widehat\phi = 11.133$ ({\bf beta regression}). The beta fit takes the influential case as indication of large dispersion.  This mistake does not occurs with ({\bf simplex regression}), 
$\widehat\sigma^2$ = 0.0346, small dispersion. As the simplex model is a dispersion model, the estimation process of the dispersion parameter  was not affected by the influence of the outlier.
However, we still  improve this simplex regression model. In this case the improvement is reached when we consider the joint modeling of the mean and the dispersion.
\begin{table}
	\caption{Parameter estimates, standard errors (s.e.), relative changes in estimates and in standard errors
		due to cases exclusions and respective $p$-values. Reading accuracy data. Simplex and Beta models with constant dispersion.}\label{tab1} 
	\renewcommand\arraystretch{1.3}
	\renewcommand{\tabcolsep}{0.17pc} 
	{\scriptsize
		\begin{tabular}{{c|c|c|c|c|c|c||c|c|c|c|c}}\hline
			Models                &\multicolumn{6}{c||}{Beta}&\multicolumn{5}{c}{Simplex}  \\ \hline  
			Parameter            &        &$\beta_1$& $\beta_2$ & $\beta_3$  & $\beta_4$ & $\phi$ &  $\beta_1$ & $\beta_2$ & $\beta_3$ & $\beta_4$ &$\sigma^2$\\  \hline
			Description &  &{Const}  &{Dislex}   & {IQ}      & {Int}            & &{Const}  &{Dislex}   & {IQ}      & {Int}            \\ \hline   
			{Full}&estimat. &1.334  &$-0.974$& 0.161 &$-0.219$&11.133&1.207  &$-0.818$&0.577 &$-0.630$&0.035\\                        
			{dataset}&e.p.     &0.138 &0.134&0.134&0.135&2.444&0.209&0.209&0.189 &0.189&0.005 \\                    
			&$p$-value&0.000&0.000&{\bf 0.238} &0.105&      &0.000&0.000&0.002&0.000&       \\ \hline              
			{Obs.1}      &change est.&$-6.0$ &$-8.4$  &65.6   &48.5    &7.7   &$-10.5$&$-15.5$&$24.4$ &22.1  &$17.7$ \\          
			{deleted}   &change s.e.&$-2.4$ &$-2.1$  &1.7    & 1.8    &8.6   &$-0.03$&$-0.03$&$-2.5$ &$-2.5$&$19.1$ \\          
			&$p$-value&0.000&0.000&{\bf 0.051} &0.018  &      &0.001 &0.000 &0.000 &0.000&       \\ \hline   
			{Obs.}&change est.&$-22.0$&$-31.3$ &204.8  &151.5   &51.8  &$-30.8$&$-45.5$&$56.2$ &50.9&99.5     \\          
			{1, 4, 7, 8}&change s.e.&$-16.3$&$-15.6$ &$-8.6$ &$-8.5$  &57.3  &$-3.4$ &$-3.4$ &$-9.5$ &$-9.5$&109.2 \\          
			{deleted}   &$p$-value&0.000&0.000 &{\bf 0.000} &0.000&      &0.000&0.027&0.000&0.000&       \\ \hline   
			{Obs.}        &change est.&12.4 &15.8&25.1&19.3&17.4        & 51.8  &  3.0   &   36.7 & 198.4&$-3.8$ \\         
			{11, 13, 16}  &change s.e.&$-8.1$&$-7.8$&$-8.4$&$-8.4$&4.4  &$-11.2$& $-4.3$ & $-8.6$ & 36.5 &3.6   	   \\    
			{deleted}    &$p$-value&0.000&0.000&0.102&0.034&      &0.000& 0.000&  0.003&0.001&                  
			\\\hline   		\end{tabular}}
\end{table}
Recentely, \cite{BAYER2017}  based on several schemes of model selection, concluded that a good model based on beta regression  for these data would have considered the following linear predictors: 
\begin{equation}\label{BMBayer:I}
\log\left(\frac{\mu_t}{1-\mu_t}\right) = \beta_1+\beta_2x_{t2}+\beta_3x_{t3}^2+\beta_4x_{t2}x_{t3}^2,
\end{equation}
and
\begin{equation}\label{BDBayer:I}
\log(\phi_t) = \gamma_1+\gamma_2x_{t2}+\gamma_3x_{t3}+\gamma_4x_{t3}^2+ \gamma_5x_{t2}x_{t3},\quad t=1,\ldots,44.
\end{equation}
We consider several simplex  and beta regression models with varying dispersion. We chose the models that present the better selection criteria values \cite{Bayer+Cribari_2015}, the better diagnostic plots and whose estimation process showed robust to the influential points. The selected models consider the following linear predictors for the mean and dispersion submodels, respectivelly,  
{\bf Simplex Model}: $g(\mu_t)=\beta_1+\beta_2x_{t2}+\beta_3x_{t3}^2+\beta_4x_{t2}x_{t3}^2$, 
$h(\sigma^2_t)=\gamma_1+\gamma_2x_{t2}+\gamma_3x_{t3}$. 
{\bf Beta Model}:
$g(\mu_t)=\beta_1+\beta_2x_{t2}+\beta_3x_{t3}^2+\beta_4x_{t2}x_{t3}^2$, $h(\phi_t)=\gamma_1+\gamma_2x_{t2}+\gamma_3x_{t3}+\gamma_4x_{t3}^2$.
We also consider that $g(\cdot) $ and $h(\cdot)$ are the logit and logarithmic  functions, respectively. 

The estimates of the parameters, standard errors and $p$-values associated with the tests of the model  parameters are, respectivelly for 
{\bf simplex distribuition}: $\widehat\beta_1$=$1.2$(0.1748) (0.0000),	$\widehat\beta_2=-0.8$(0.1748) (0.0000), $\widehat\beta_3=0.4$(0.0739)(0.0000),
$\widehat\beta_4$=$-0.4$(0.0739) (0.0000), $\widehat\gamma_1=1.1$(0.2162 )(0.0000), $\widehat\gamma_2=-2.8$(0.2633) (0.0000), 
$\widehat\gamma_3=-0.6$(0.2639) (0.0222). For 
{\bf beta distribuition}: $\widehat\beta_1$=$1.1$(0.1509)(0.0000),	$\widehat\beta_2=-0.8$(0.1509)(0.0000), $\widehat\beta_3=0.4$(0.0699)(0.0000),
$\widehat\beta_4$=$-0.4$(0.0698)(0.0000),
$\widehat\gamma_1=1.1$(2.6)(0.0000), $\widehat\gamma_2=1.2$ (0.2623) (0.0000), 
$\widehat\gamma_3=1.0$(0.2508)(0.0000), and $\widehat\gamma_4=0.9(0.2092)(0.0000)$.

\begin{figure}
	\centering
	\includegraphics[width=0.52\linewidth, angle=270]{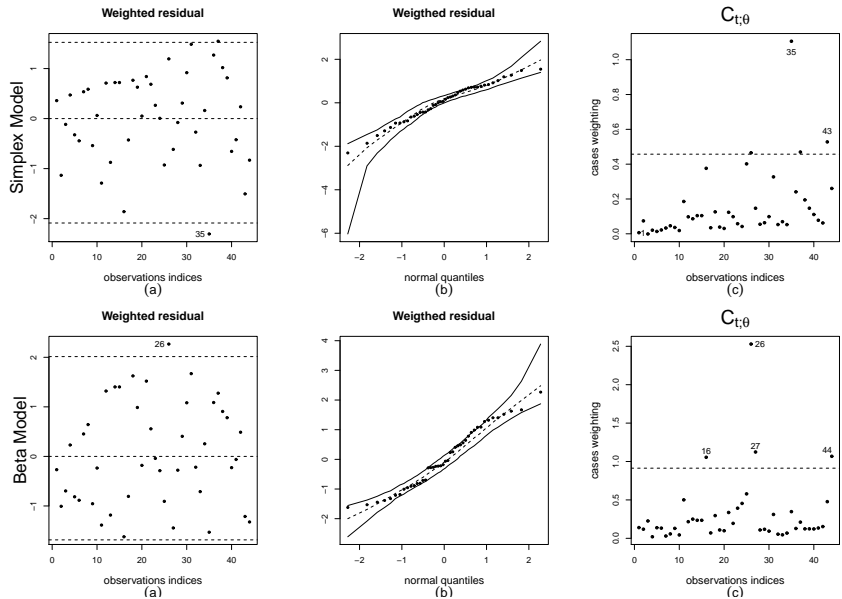}
	\caption{Diagnostic plots. Reading accuracy data. Simplex model: $\log({\mu_t}/{1-\mu_t}) = \beta_1+\beta_2x_{t2}+\beta_3x_{t3}^2+\beta_4x_{t2}x_{t3}^2$ and 
		$\log(\sigma^2_t) = \gamma_1+\gamma_2x_{t2}+\gamma_3x_{t3}$.  Beta model: $\log({\mu_t}/{1-\mu_t}) = \beta_1+\beta_2x_{t2}+\beta_3x_{t3}^2+\beta_4x_{t2}x_{t3}^2$ and 
		$\log(\phi_t) = \gamma_1+\gamma_2x_{t2}+\gamma_3x_{t3}+\gamma_4x_{t3}^2$. }
	\label{fig:Hetero}
\end{figure}

We noticed that all covariates were statistically significant, both for the mean model and for the dispersion model. The beta regression model that presents the better diagnostic plots is similar to the model presented in \ref{BMBayer:I} and \ref{BDBayer:I}.  However, it does not consider the interation, on dispersion submodel. 
In Figure~\ref{fig:Hetero} we present the diagnostic plots for the final candidate models.  

This figure reveals that the simplex distribution seems to fit very well to this data since the most of all points are randomly distributed within  the envelope bands 
(Figure~\ref{fig:Hetero}(b)). 
Otherwise, the normal probability plot of the beta regression presents a slight trend of the positive residuals  to being out of the upper band of the simulated  envelope (Figure~\ref{fig:Hetero}(e)).
The Figures~\ref{fig:Hetero}(c) and ~\ref{fig:Hetero}(f) show that  case 26 can be more influential for the estimation process of beta regression than  case  35  can be influential for the simplex model. In fact, the maximum normal curvature $C_{\rm max}$ (cases weighting) for the simplex and beta regression models are, respectively, 1.1 and 2.0. To access influence in beta regression we use the measures developed by \cite{ESPINHEIRA2008Inf}  and \cite{Ferrari+Espinheira+Cribari_2011}.

The relationships between scores in reading accuracy, status of dyslexia and non verbal IQ are used as protocol to make several diagnostics  in the psychological area \cite{SMITHSON2006}.  Thus,  we consider very important make suggestions. For while, the simplex regression model presented in Figure~\ref{fig:Hetero}  is a better option for modeling the reading accuracy data than the beta regression model. 
Our conclusion is based on the fact that the simplex distribuition fits the data better than the beta distribuition (Figures~\ref{fig:Hetero}(b) and ~\ref{fig:Hetero}(e)) and because the likelihood  estimation scheme  is more robust for influential cases when we consider the simplex regression than when we consider the beta regression.
\section{Application II : Fluid Catalytic Cracking (FCC) data.}
The data of this application come from a work to obtain the degree in chemistry from the Faculty of Sciences, Department of Chemistry of the National University of Colombia, Salazar's graduation work in 2005. The FCC (Fluid Catalytic Cracking) process is used to convert hydrocarbons of high molecular weight into small molecules of higher commercial value, through the contact of these with a catalyst.
The FCC process is often considered the heart of a refinery since it allows the production to be adapted to the products of greater demand and  high profitability. The main catalyst of the process  is zeolite USY. 
Other important substance that participates in the catalyzing process is vanadium. Vanadium is known to participate in the destruction of the catalyst by reducing the active surface, the selectivity and the crystallinity of zeolite USY especially in the presence of steam. 	It is known that every 1000 ppm of vanadium on the catalyst decreases gasoline yields by about $2.3\%$ \citep{Salazar2005}. 
The process also depends on the temperature, which should be close to 720 degrees celsius. 
The interest is in  modelling the percentage of crystallinity of the zeolite USY $(y)$ based on diferent concentrations of vanadium $(x_2)$ and of  steam $(x_3)$, considering two values for the temperature of the process $(x_4)$. It is expected that the larger the concentrations of vanadio and steam, the smaller will be the crystallinity percentage.

At first  we consider linear simplex models with varying dispersion and the final candidate that includes all covariates that were statistically significant was
\begin{equation}\label{PCR:Linear_1}
\log\frac{\mu_t}{1-\mu_t}=\beta_1+\beta_2x_{t2}+\beta_3x_{t3}+\beta_4x_{t4};\quad\log{\sigma^2_t}=\gamma_1+\gamma_4x_{t4},\,t=1,\ldots, 28.
\end{equation}
We also consider a beta regression model with the same predictors and link functions presented in \eqref{PCR:Linear_1}. In Figure~\ref{fig:Linear} we presented the residual plots against the observation indexes. The nonlinear trend of the residuals becomes evident both in the simplex model and in the beta model (Figures~\ref{fig:Linear} (a) and (c)).
\begin{figure}[!ht]
	\centering
	\includegraphics[width=0.52\linewidth, angle=270]{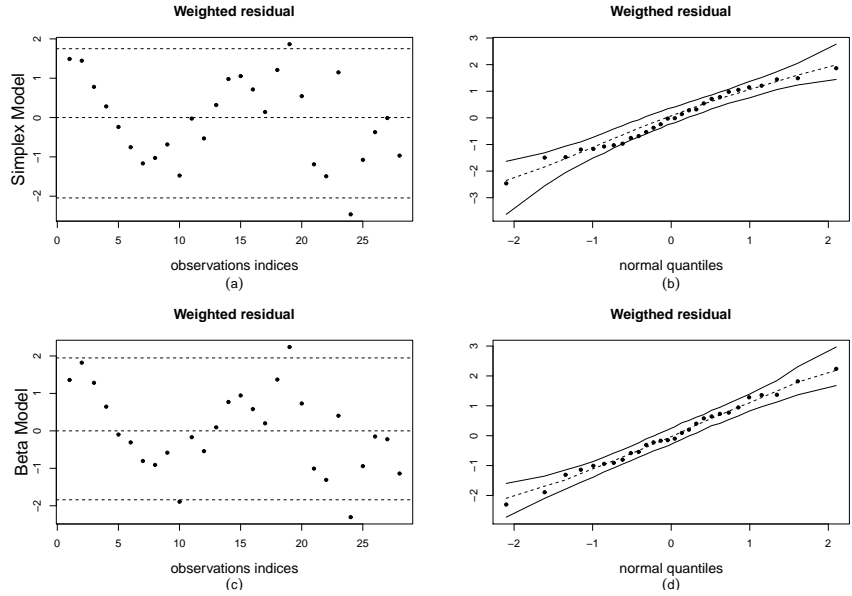}
	\caption{Diagnostic plots. Linear Models (FCC) data.  }
\end{figure}\label{fig:Linear}
This data is being modeled for the first time. There are no suggestions of non-linear functions
to this problem. Thus, initially we present a descriptive analysis of the response and covariates using boxplots (Figure~\ref{fig:Boxplot}). 		
\begin{figure}
	\centering
	\vspace{-2.02cm}
	\includegraphics[width=0.62\linewidth, angle=270]{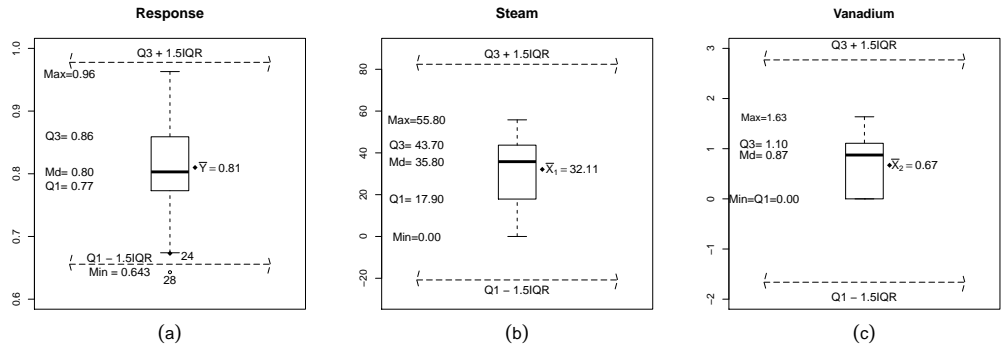}
	\vspace{-2.1cm}
	\caption{Boxplots. FCC data.}
\end{figure}\label{fig:Boxplot}
Based on Figure~\ref{fig:Boxplot}(a) it can be noted that the observed response is concentrated near the upper end of the standard unit interval. In fact, $75\%$ of observations are equal to or greater than 0.77. Additionally, {\bf the presence of an outlier is observed}. Observation 28 assumes the lowest response value, which is equal to 0.643 near to value of case 24, which is 0.674.  In Figure~\ref{fig:Boxplot}
(b) it is noted a high dispersion for the variable steam, with occurrence of extreme values, such as zeros, which is associated to the fact that $50\%$ of the observations are greater than 35.8. It is worth noting that the variable could assume plausible values between $-20$ and 80. Considering a regression model including steam as a covariate, it seems interesting to use a function in order to stabilize this dispersion. A plausible function is the logarithm.
However, the presence of zeros implies that we would have to consider, for example, $log(\text Steam + \beta)$, i.e, a nonlinear function in the predictor. Another possibility is to consider function such as $1/(\text Steam + \beta)$, $\text Steam/(\text Steam + \beta)$, which are also nonlinear functions. Finally, including vanadium as an additional covariate, we could think of the same functions used to establish the dispersion of steam, but since the first quartile is equal to the minimum and equals to zero for vanadium, a plausible function would be $\sqrt{\text {Vanadium}}$.

Based on this discussion, we evaluated some competing models, both based on the simplex distribution and on the beta distribution, where the models that presented the best residuals plots were those that considered the following predictors:
\begin{equation}\label{PCR:NoLinear}
\eta_t=\beta_1+\beta_2\frac{x_{t2}}{x_{t2}+\beta_3}+\beta_3x_{t3}+\beta_4\sqrt{x_{t4}};\quad\zeta_t=\gamma_1+\gamma_2x_{t4}^2,\, t=1,\ldots,28.
\end{equation}
The next step is the construction of the inicial guess. Based on \eqref{PCR:NoLinear} we have that the $th$ row the matrix $\tilde{\mathcal{X}}^{(0)}=[\partial \eta/\partial \beta]_{\beta = \beta^{(0)}}$ is defined as
$(1,x_{t2}/(x_{t2}+\beta_{3}^{(0)}),-\beta_{2}^{(0)}x_{t2}/((x_{t2}+\beta_{3}^{(0)})^{-2},x_{t3},x_{t2}^2)$, $t=1,\ldots, 28$. 
Note that in \eqref{PCR:NoLinear} we have more parameters than covariates. To obtain $(\beta_{1}^{(0)},\beta_{2}^{(0)},\beta_{4}^{(0)}, \beta_{5}^{(0)}),^{\top}$ we could use $(X^{\top}X)^{-1}X^{\top}g(y)$, in which the $t$th row of matrix $X$ can be defined as  $x^{\top}_t=(1,x_{t2}/(x_{t2}+\beta_{3}^{(0)}),x_{t3},x_{t2}^2)$, $t=1,\ldots, 28$. However, it is necessary to provide a numerical value for $\beta_{3}^{(0)}$. 
Based on the boxplot (Figura~\ref{fig:Boxplot}) and in the nonlinear functions, we will assign plausible values for $\beta_{3}^{(0)}$. Consider the equations $[{\text {steam}} + \beta_{3}^{(0)} = -20]$ and $[{\text {steam}} + \beta_{3(0)} = 80]$. Taking steam$=0$, it follows that 
$\beta_{3}^{(0)}\in (-20, 80)$. Taking steam$=55.80$, it follows that 
$\beta_{3}^{(0)}\in (-75.80, 24.2)$. 
Chosen possible values for $\beta_{3}^{(0)}$ we obtain
$\beta_{NL}^{(0)} = (\tilde{\mathcal{X}}^{(0)\top}\tilde{\mathcal{X}}^{(0)})^{-1}\tilde{\mathcal{X}}^{(0)\top}(g(y)-f(x,\beta_{L}^{(0)}))$, where  
$f(x_t^{\top},\beta_{L}^{(0)}) = \beta_{1}^{(0)} + \beta_{2}^{(0)}x_{t2}/(x_{t2}+\beta_{3}^{(0)})+ \beta_{4}^{(0)}x_{t3} + \beta_{5}^{(0)}x_{t2}^2$. 
We tested same values $\beta_{3}^{(0)}$ and we chose $\beta_{3}^{(0)}=-20$  that led to parameter estimates with finite standard errors for both the simplex and the beta models. After the convergence of the process we obtained the model estimates presented in Table~\ref{tab3}. It is important to emphasize that the construction of diagnostic plots is only possible after the definition of the starting values procedure as described above. Thus, to select the best models based on the residual analysis, we carried out the complete estimation process considering starting values for all models that we consider as candidates. The chosen model is defined in  \eqref{PCR:NoLinear}.

\begin{figure}[!ht]
	\centering
	\includegraphics[width=0.52\linewidth, angle=270]{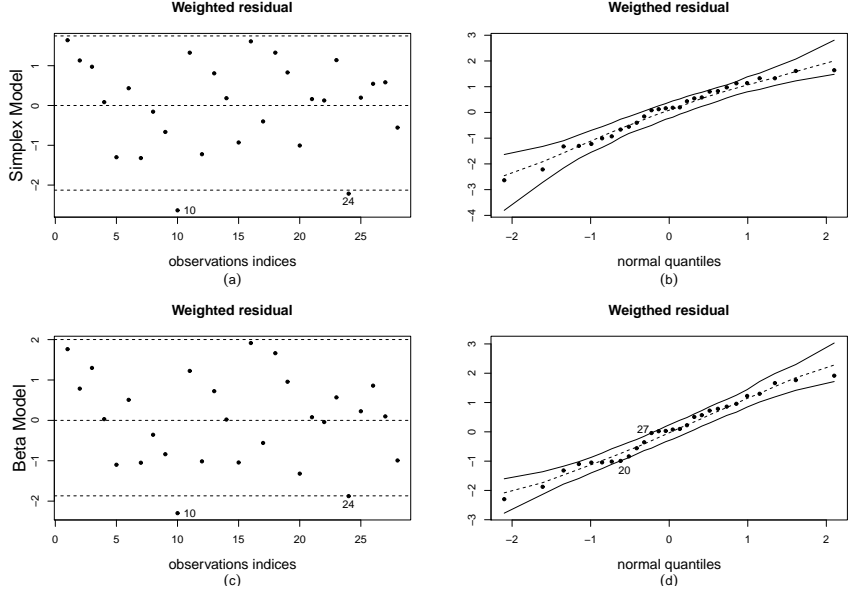}
	\caption{Diagnostic plots.  Simplex and Beta models. (FCC) data.}\label{fig:Hetero1}
\end{figure}

\begin{figure}[!ht]
	\centering
	\includegraphics[width=0.52\linewidth, angle=270]{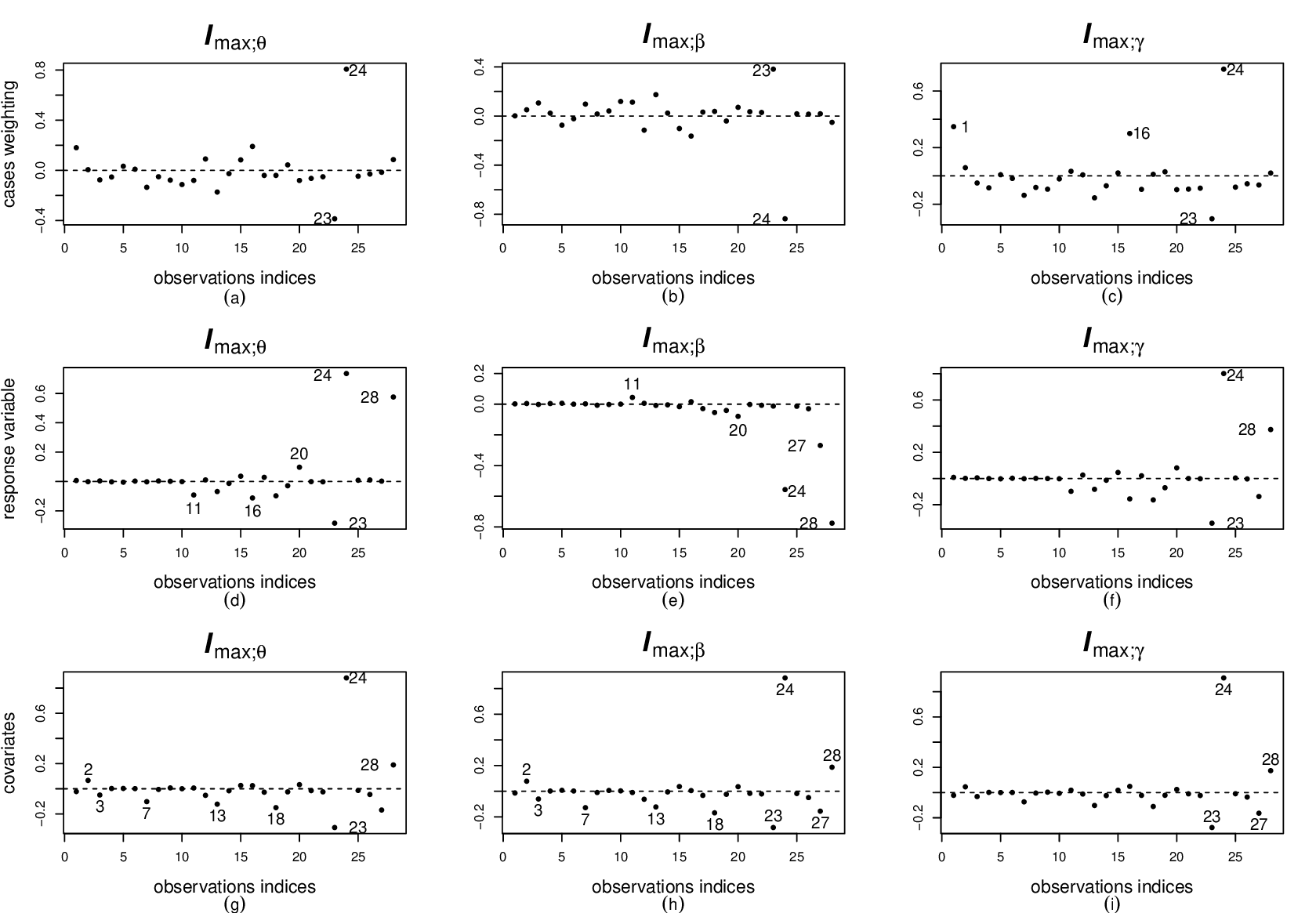}
	\caption{Local influence plots $I_{\text{max}}$. Simplex Model. (FCC) data.  }\label{fig:Hetero2}
\end{figure}

\begin{figure}[!ht]
	\centering
	\includegraphics[width=0.52\linewidth, angle=270]{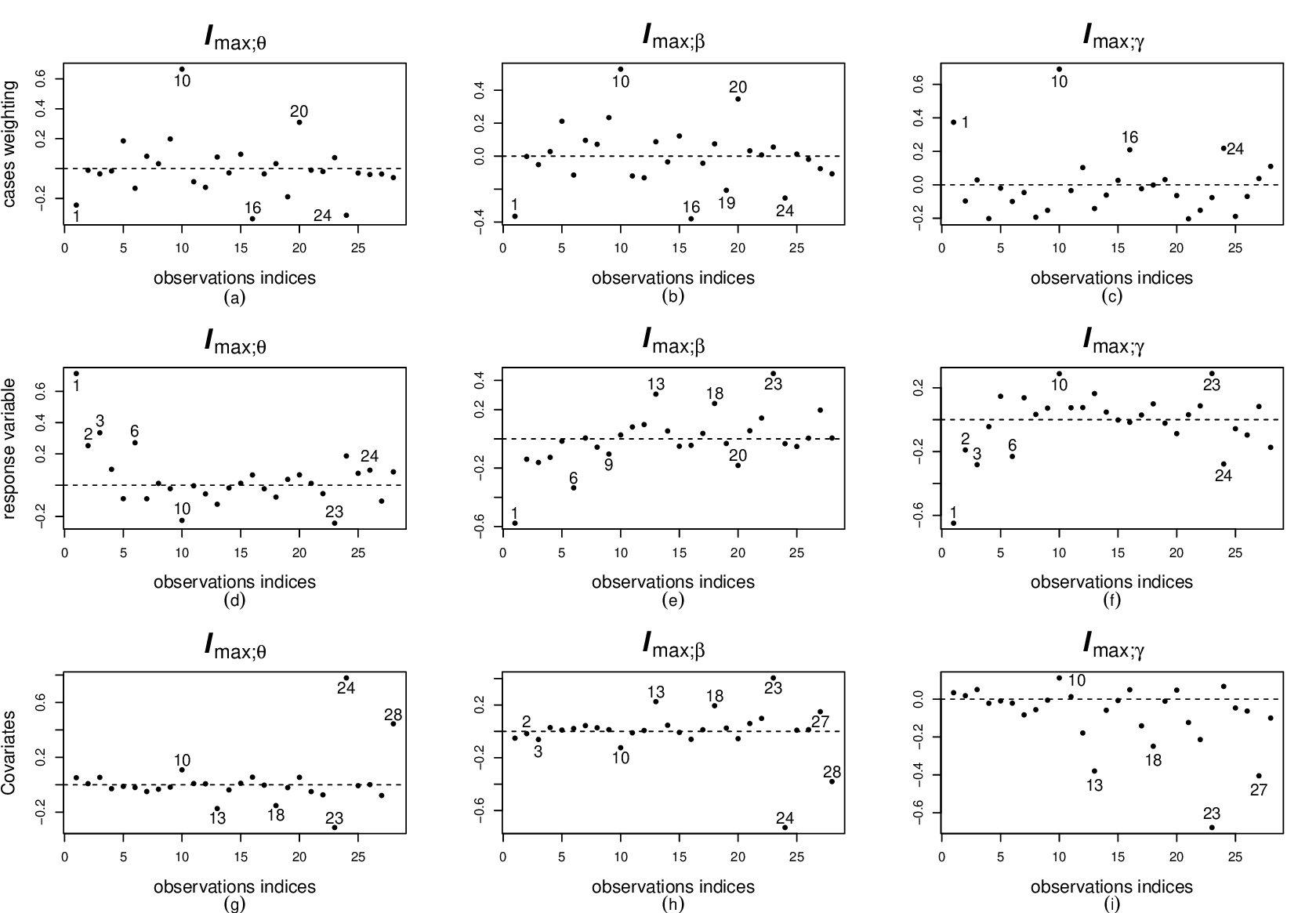}
	\caption{Local influence plots $I_{\text{max}}$. Beta Model. (FCC) data.  }\label{fig:Hetero3}
\end{figure}

\begin{figure}[!ht]
	\centering
	\includegraphics[width=0.52\linewidth, angle=270]{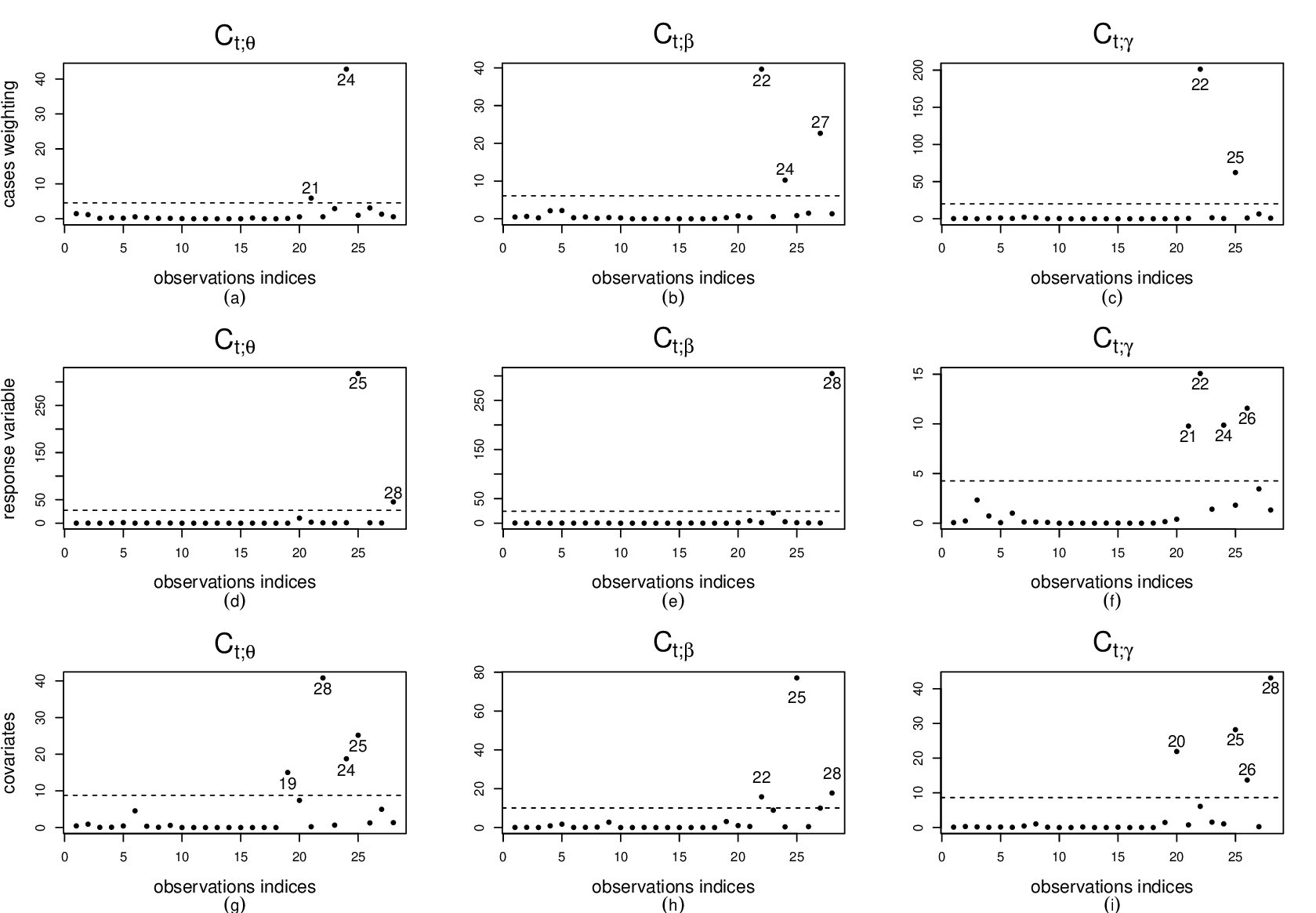}
	\caption{Local influence plots, $C_t$. Simplex Model. (FCC) data.  }\label{fig:Hetero4}
\end{figure}

\begin{figure}[!ht]
	\centering
	\includegraphics[width=0.52\linewidth, angle=270]{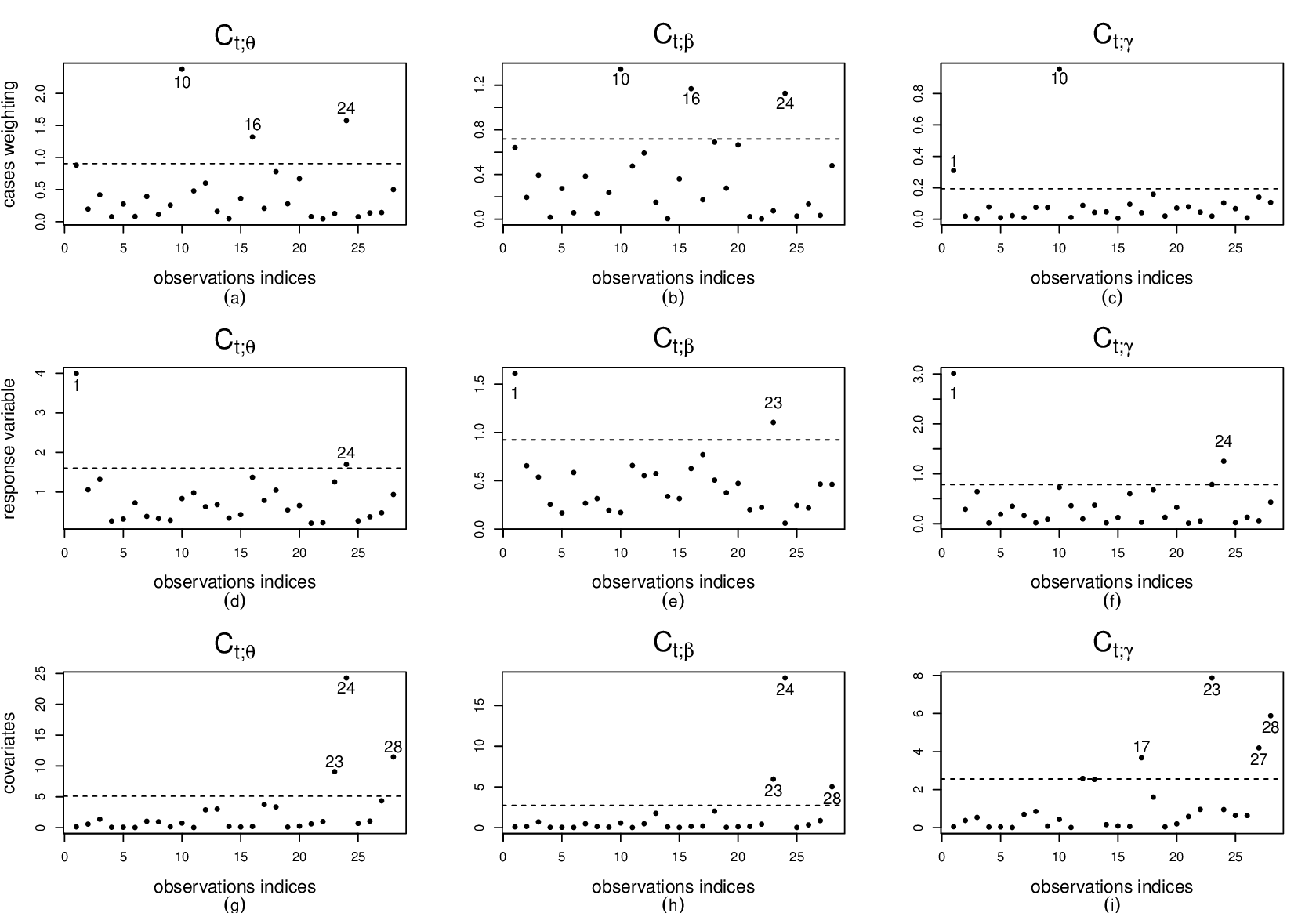}
	\caption{Local influence plots, $C_t$. Beta Model. (FCC) data.  }\label{fig:Hetero5}
\end{figure}

After the selection of the simplex and beta models throughout residual analysis,  we carried out the influence analysis to verify the robustness of the maximum likelihood estimation procedure for each model in the presence of influential points. We have to stress that, to carry out the influence analysis on nonlinear beta regression models, we use the results presented by \cite{ROCHA2011}. In Figure~\ref{fig:Hetero1}, we present the residual plots for the simplex and beta models. We can verify that both models fit the data well and that the nonlinear trends that appeared in the residual plots against the indices of the observations is considerably smoothed after the inclusion of nonlinear predictors, as shown in Figures~\ref{fig:Hetero1} (a) and (b). In addition, the normal probability plots with simulated envelopes show that the adherence of the data to the simplex distribution is slightly better than that observed with using the beta distribution. We also highlight there are two points that lie within the boundaries of the envelopes bands (cases 20 and 27) in the case of the beta regression. Regarding the outliers,  observations 10 and 24 stood out for both regressions.
In Figures~\ref{fig:Hetero2} through~\ref{fig:Hetero5}, we present the local influence and total local influence plots for the simplex and beta regression models. We consider the case weighting (first line); response perturbation (second line) and simultaneous disturbance of the vanadium covariate in both mean and dispersion models (third line). The influence plots identify important cases for the vector $\theta=(\beta^{\top},\gamma^{\top})^{\top}$ (first column), for the vector $\beta^{\top}$ in particular (second column), and for the $\gamma^{\top}$(thrid column).

In Figures~\ref{fig:Hetero2} and \ref{fig:Hetero3}, the measures  $I_{max}$ reveal joint influential observations, for simplex and beta models, respectively. We highlight that the cases with opposite signs in  $I_{max}$ generally exert opposite influences on the estimation of the regression models. This fact indicates that these cases should be studied jointly. For example, in general, the observations 23 and 24 should be investigated separately (Figura~\ref{fig:Hetero2})(a)(simplex model).  Similarly, for the beta regression, the sets $\{1,16,24\}$ and $\{10,20\}$ must be evaluated separately. To investigate the real influence of theses cases, they are excluded from the data and the model is reestimated. In general, several sets should be evaluated, among which we highlight 
$\{1, 16, 24\}$, $\{20,24,28\}$, $\{20,24,27,28\}$, $\{13, 18, 23, 27\}$ , $\{24, 28\}$. 
It is also interesting to evaluate cases that individually exert a disproportional influence on the model fit. This information can be assessed by analysis of plots of $C_t$, (Figure~\ref{fig:Hetero4}). Several points should be investigated individually, including specially cases 22, 24, 25, 27 and 28, for the simplex regression.
For the beta model, the influential cases 1, 10, 16, 17, 24, 27 emerge individually ( Figure~\ref{fig:Hetero5}). In Table~\ref{tab3} we present the changes on the  parameter estimates, in the estimates of standard errors and the p--values related to the significance tests of the parameters after the exclusion of cases and sets of cases that were, in fact, more influential. 
\begin{table}
	\caption{Parameter estimates, standard errors (s.e.), relative changes in estimates and in standard errors
		due to cases exclusions and respective $p$-values. Mean Model: $\beta_1 + \beta_2(\text{Steam} + \beta_3)^{-1} + \beta_4\text{InTemp} + \beta_5\sqrt{\text{Vanadium}}$ and Dispersion Model: $\gamma_1+\gamma_2\text{Vanadium}^2$.}\label{tab3}
	\renewcommand\arraystretch{1.3}
	\renewcommand{\tabcolsep}{0.0000005pc} 
	{\scriptsize
		\begin{tabular}{{c|c|c|c|c|c|c|c||c|c|c|c|c|c|c}}\hline
			{Model}&\multicolumn{7}{c||}{Simplex}&\multicolumn{7}{c}{Beta}\\\hline 
			&$\beta_1$& $\beta_2$ & $\beta_3$  & $\beta_4$ & $\beta_5$ &  $\gamma_1$ & $\gamma_2$&$\beta_1$& $\beta_2$ & $\beta_3$  & $\beta_4$ & $\beta_5$ &  $\gamma_1$ & $\gamma_2$ \\  \cline{1-15}
			{Full} & 2.374  & $-0.106$& $-27.8$& $-0.290$ &  $-0.751$& 0.8246 & $-1.217$&2.296 &$-0.101$&$-26.9$&$-0.286$&$-0.679$&3.915&0.400\\        	                       
			{data}& 0.149  & 0.042   & 4.021    & 0.106    &  0.141   & 0.369  &  0.314	 &0.152 & 0.039 & 4.358 & 0.107 & 0.129 &0.369 &0.310	  \\             
			&0.000&0.012& 0.000&0.006&0.000&0.025&0.000&0.000&0.010&0.000&0.007&0.000&0.000&0.200\\\hline          
			\multicolumn{15}{c}{Relative Changes $(\%)$ in estimates and in standard errors due exclusions. $p$--values after exclusions. }\\\hline
			{Obs.}	 &$-4.4$&   0.2 &$-3.8$&  3.1 &$-10.3$&$-55.6$&$-26.5$	&$-4.0$ &$-3.6$ & $-3.2$&$-2.6$&$-9.1$& 5.4  &$-41.4$   \\     
			{1}      &$-0.8$&$-20.7$& 2.01 &  2.6 &$-10.3$&  3.6  &  1.7	& 1.5   &$-15.7$& $4.5$ &$0.6$ &$-5.7$& 3.5  & 1.7		 \\
			{deleted} &0.000&0.001&0.000&0.006&0.000&0.339&0.00&0.000&0.003&0.000&0.010&0.000&0.000&0.461\\\hline
			{Obs.}	   &3.6&$-9.2$&2.4&$-7.9$&10.8  &$-19.6$& $-4.5$&4.2  &$-6.9$&2.3   &$-5.9$& 11.9 & 10.6 & $-52.0$  \\                     
			{10}       &$-7.2$&2.8&$-2.1$&$-3.1$&$-6.0$& 3.6   & 1.7 &$-9.6$ &$-0.5$&$-4.5$&$-5.3$&$-8.7$&3.7   &  1.7  \\                    
			{deleted} &0.000&0.026&0.000&0.010&0.000&0.084&0.000&0.000&0.017&0.000&0.008&0.000&0.000&0.547 \\\hline            
			{Obs.}	  &$-1.0$&$-8.6$&$-1.1$&$-6.7$&$ -0.4$&$-82.9$&$-34.0$& $-0.2$&$-10.3$&$-1.1$&  $-9.3$&  $ 1.7 $& $18.3$&$-107.5$  \\ 
			{1,10}    &$-10.3$&$-16.9$&$ 0.2$&$-1.9$&$-17.0$&$  7.6$&$  3.7$&$-10.5$&$-18.4$&$ 0.3$&  $-7.2$&  $-15.8$& $7.6 $&$  3.6 $   \\
			{deleted} &0.000&0.006&0.000&0.010&0.000&0.724& 0.014&0.000&0.005&0.000&0.009&0.000&0.000&0.928\\\hline
			{Obs.}	  &0.9&$-12.8  $ &  $-1.3$ & 8.8  &$-4.6$ &$ 3.9  $&$21.3  $&$ 0.1  $&$-6.3  $&$-0.6$&$  25.0$&$-7.4$& $-0.4$&42.2 \\  
			{24}      &$-3.0  $&$-15.5  $ & $ 7.7$ & $-5.4$& $-1.1$&$ 1.4  $&$10.2  $&$-3.3  $&$-8.8  $&$  1.3  $&$-1.6$&$-0.7$& $1.3$&10.2 \\  
			{deleted} &0.000&0.009&0.000&0.000&0.000&0.022&0.000&0.000&0.009&0.000&0.001&0.000&0.000&0.099 \\\hline  
			{Obs.}	  & $-3.0$&$-21.7$&$-2.9$&$ 9.8$&$-12.8$&13.5  &37.5	   & $-3.1$ &$-14.0$&$-1.5$ &$ 4.6$ &$-15.0$&$-2.5$&93.2  \\  
			{24,28}   & $ 6.0$&$-23.5$&$16.5$&$-5.7$&$  3.4$& 3.7  &25.9	   & $ 2.0$ &$-16.8$&$ 6.7$ &$-4.1$ &$  3.1$&$ 3.7$&25.9   \\  
			{deleted} &0.000&0.010&0.000&0.002&0.000&0.015&0.000&0.000&0.008&0.000&0.004&0.000&0.000&0.050\\\hline  
			{Obs.}	  &$-3.3$&$-26.6$&  2.0 & $-5.8$& $-16.6$& 16.9 & 47.1	 &$-2.9$ &$-20.2$&  4.4  &$-13.3$& $-17.7$& $-3.0$&141.6   \\  
			{20,24}   &$ 3.8$&$-9.1$& 12.6 & $-9.4$& $4.3$& 3.8  & 27.5	 &$-3.1$ &-$1.7$&  2.1  &$-10.6$& $  2.7$& $ 3.7$& 27.6 	   \\  
			{28,del.} &0.000&0.042&0.000&0.005&0.000&0.012& 0.000&0.000&0.038&0.000&0.010&0.000&0.000&0.015\\\hline  
			{Obs.}	  &$-1.0$&$0.4$ &$  6.0$&$-39.6$&1.2   &44.2  &47.5  &$-0.1$&5.1&$2.0 $ &$-19.6$ &$4.6$ &$-3.0$&17.2	     \\  
			{13,23}   &$8.4$&$6.5$ &$-27.1$&$-19.3$&10.4  &0.7   & 11.5 &$5.6$ &10.6  &$-4.4$ & $ 12.5$&$6.4$ &$ 0.7$&11.5	     \\     
			{27,del.} &0.000&0.0174&0.000&0.042&0.000&0.001&0.000&0.000&0.015& 0.000& {\bf 0.058}&0.000&0.000&0.179     \\\hline
			{Obs.}	      &$-1.3$ &$-18.1$&12.3   &$-54.3$&$-3.9$&60.3 &71.3 &1.7  &$-7.1$      &11.2   &$-54.6$&9.1  &$-8.7$&121.8      \\  
			{13,20,23}       & 12.8  &5.79   &$-39.6$&$-35.6$& 15.1 &0.9  &12.6 &  5.3&   26.5     &$-23.9$& -4.3  &10.3 &10.0  &12.7    \\     
			{27,del.}  &  0.000&0.050  &0.000  &0.053  & 0.000&0.000&0.000&0.000& {\bf 0.061}&0.000  & {\bf 0.210} &0.000&0.000 &0.012     \\\hline
	\end{tabular}}                                 
\end{table}

All the parameters are statistically significant, both for mean and for the dispertion submodels, in the simplex regression (Table~\ref{tab3}). Such conclusion is not corroborated by using the beta regression because the coefficient of covariate  $\text{Vanadio}^2$ covariate is not statistically significant ($p$--value$=$0.200). This result is even more compeling due to the exclusions of the case 1 (p--valor$=$0.461), of the case 10; (p--value$=$0.557) and of the set
$\{1,10\}$; ($p$--value$=$0.927). On the other hand, the unique exclusion of case 24 goes on the opposite direction, to the point of implying on the significance of $\text{Vanadio}^2$, to the value of $10\%$, ($p$--value$=$0.099). This fact is marked due to the exclusion of cases $\{24,28\}$, ($p$--value$=$0.0497), $\{20,24,28\}$ ($p$--value$=$0.0154) and $\{13,20,23,27\}$ ($p$--value$=$0.0120). 

The observations $\{13,23,27\}$ interfere on the estimation of the mean submodel with greater impact for {\bf the beta regression}. The exclusion of the set $\{13,23,27\}$ implies on $p$--values associated to the covariate coefficient related to temperature, which are equal to 0.0420, for the simplex model and   0.0580,  for the beta model. This fact would entail on the nonsignificance of this covariate, at $5\%$ level, on the mean submodel of the beta regression, while it is statistically significant, at $1\%$ level, when considering the complete data. This result on {\bf beta regression} is substantially more pronounced when the set $\{13,20,23,27\}$ is excluded from the data, $p$--value$=$0.210 for $\beta_4$ (temperature) and $p$--value$=$0.061 for $\beta_2$ (steam). The inferential results of the simplex model  are also affected by the set $\{13,20,23,27\}$. However, both  coefficients of steam and of temperature are near of $5\%$ level of significance.

Once again  the dispersion of the beta model is affected by outliers and influential points. For the complete data, we have that ${\widehat \phi}_{\text{max}} = 147$. However, when the cases 24 and 28 are excluded of the data, the dispersion decreases considerably, such that maximum ${\widehat \phi}_{\text{max}} = 235$. The same does not occur with the simplex model, the maximum dispersion with or without the influential cases is roughly the same, ${\widehat \sigma^2}_{\text{max}}\approx 2.4$ .
\section {Concluding remarks}\label{S:conclusions} 
On this paper we propose a simplex regression model which considers nonlinear structures for the parameters for both the mean and the dispersion submodels. We propose a technique for defining starting values for Fisher's iterative scoring process for the nonlinear simplex regression model
Additionally, we propose a residual based on Fisher's iterative process to estimate the parameter vector $\beta$and local influence measures considering: case weighting, disturbance of the response variable and the joint disturbance of countinuous covariates for both the mean and the dispersion submodels.

We present Monte Carlo simulations to investigate the empirical distribution of the proposed residual and we verified a good approximation of this distribution by the standard normal distribution. However, we highlight that to define the residual plot against the observation index, we used the proposal defined by \cite{ESPINHEIRA2017} based on empirical residual quantiles obtained from the simulated envelope  of the standard normal distribution.

We also presented two applications to real data: one linear and the other nonlinear on the parameters. On these applications we compared the fit of the simplex and  beta regression models. On both applications the simplex regression model showed to be a better option than the beta regression model. In  special, because the estimation process by maximum likelihood associated with the simplex model showed to be more robust than the beta model,  when there are highly influential cases in the data.

\end{document}